\documentstyle[amssymb,amsmath,eucal]{article}

\begin{document}

%\righthyphenmin=2
%\textwidth=180mm
%\textheight=240mm

\def\1{{ }}

\def\s{\scriptstyle\,\,\,}
\def\t{\textstyle}
\def\scc{\sc {}}
\def\bff{\bf {}}
\def\bfff{\large\bf }

\def\half{{\t \frac12}}

\renewcommand\S{Section }
\def\graph{{\rm GRAPH}}

\def\R {{\Bbb R }}
\def\C {{\Bbb C }}
\def\K{{\Bbb K}}
\def\H{{\Bbb H}}

\def\U{{\rm U}}

\def\const{{\rm const}}
\def\B{{\rm B}}
\def\Gr{{\rm Gr}}
\def\M{{\rm Mat}}
\def\r{{\frak R}}
\def\W{{\rm W}}
\def\SW{{\rm SW}}

\def\O{{\rm O}}
\def\SO{{\rm SO}}
\def\GL{{\rm GL}}
\def\SL{{\rm SL}}
\def\SU{{\rm SU}}
\def\Sp{{\rm Sp}}
\def\SOS{\SO^*}

\def\ch{\cosh }
\def\sh{\sinh}
\def\tg{\tan}
\def\th{\tanh}

\def\ov{\overline}
\def\phi{\varphi}
\def\epsilon{\varepsilon}
\def\kappa{\varkappa}%don't change to $\kappa$!!!!!!
\def\le{\leqslant}
\def\ge{\geqslant}

\newcommand{\Res}{\mathop{\rm Res}\nolimits}
\newcommand{\supp}{\mathop{\rm supp}\nolimits}

\def\Cay{{\rm Cay}}

\def\bol{\vspace{22pt}}
\def\mal{\smallskip}

\renewcommand{\Re}{\mathop{\rm Re}\nolimits}
\renewcommand{\Im}{\mathop{\rm Im}\nolimits}

\def\F{{\frak F}}
\def\L{{\cal L}}
\def\b{{\cal B}}

\def\risodin{
\begin{picture}(100,90)
\unitlength=1mm
\put(50,20){\circle*{1}}
\put(53,23){$(p+q)/2-1$}
\put(60,40){$\Pi_0$}
\put(43,40){$\Pi_1$}
\put(33,40){$\Pi_2$}
\put(20,35){\circle*{1}}
\put(20,35){\circle{4}}
\put(13,35){${\cal O}_\delta$}
\put(25,40){\vector(-1,-1){4.5}}
\put(23,41){$\alpha_0$}
%\thinlines
\linethickness{0.05mm}
\multiput(20,10)(10,0){4}
{\line(0,1){60}}
\linethickness{0.3mm}
\put(0,20){\vector(1,0){100}}
%\put(0,0){ Picture 1. The complex plane $\alpha$.}
\end{picture}
}

\def\risdva{
\unitlength=1mm
\begin{picture}(100,150)(-50,-80)
{\linethickness{0.005mm}
\put(0,-60){\vector(0,1){120}}
\put(-50,0){\vector(1,0){100}} }
%\linethickness{0.4mm}
\put(0,-20){\line(0,1){40}}
\put(0,40){\line(0,1){15}}
\put(0,-40){\line(0,-1){20}}
\put(0,30){\oval(20,20)[r]}
\put(0,-30){\oval(20,20)[l]}
\put(-1,30){\line(1,0){2}}
\put(-1,-30){\line(1,0){2}}
\put(5,30){\circle{1}}
\put(-5,-30){\circle{1}}
\put(-10,30){\circle*{1} }
\put(10,-30){\circle*{1}}
\put(-8,30){\vector(1,0){5}}
\put(8,-30){\vector(-1,0){5}}
\put(15,30){$S_1$}
\put(5,50){$R_1$}
\put(5,5){$Q$}
\put(-15,-30){$S_2$}
\put(5,-50){$R_2$}
\end{picture}}

\vspace{22pt}

\begin{center}

{\large\bf Plancherel formula for Berezin deformation
 of $L^2$ on Riemannian symmetric space}

\vspace{22pt}

Yurii A. Neretin\footnote%
{supported by grants RFBR-98-01-00303 and RFBR 96-01-96249}

\end{center}

\vspace{22pt}

{\bfff 0.Introduction}

\bol

\newcounter{sec}
 \renewcommand{\theequation}{\arabic{sec}.\arabic{equation}}

\newcounter{fact}
\def\fact{\addtocounter{fact}{1}{\scc \arabic{fact}}}

\newcounter{punkt}
\def\punkt{\addtocounter{punkt}{1}{\bff \arabic{punkt}}}

{\bff 0.1. Kernel representations.}
Let $G$ be a {\it classical} real group and let
$K$ be its maximal compact subgroup.
 Consider the Riemannian noncompact symmetric space
$G/K$. There exists a hermitian symmetric
space
$$\widetilde G/\widetilde K\supset G/K$$
such that
$$\dim_{\R} G/K =\frac12\dim_\R\widetilde G/\widetilde K$$
and $G/K$ is a totally real submanifold in $\widetilde G/\widetilde K$
(the list of embeddings $G/K\to \widetilde G/\widetilde K$
  see in Section 6).
We say that  the symmetric space
$\widetilde G/\widetilde K\supset G/K$
is the {\it hermitization} of the symmetric space
$G/K$.

We define a {\it kernel representation} $\rho$ of the group $G$
as a restriction of an unitary
 highest weight representation
 $\widetilde\rho$ of $\widetilde G$ to
the subgroup $G$. By well-known Harish-Chandra
construction, highest weight representations
of $\widetilde G$ are natural representations in  spaces
of (scalar-valued or vector-valued) holomorphic functions on
 $\widetilde G/\widetilde K$. We say that $\rho$ is
a  {\it scalar valued
kernel representation} if $\widetilde\rho$
is realized in scalar-valued
holomorphic functions on $\widetilde G/\widetilde K$.

The kernel representations are  deformations of $L^2(G/K)$
in some precise sense explained in Subsection 1.13.

The purpose of this paper\footnote%
{This work is  continuation of works
\cite{Ner4}, \cite{Ner5} but logically it is independent
on these papers.  Our main result was announced in
\cite{Ner5}.}
 is to obtain the Plancherel formula for
scalar-valued kernel
representations
(see formula (2.6)--(2.15)).

There were different reasons for interest, which was
attracted by kernel representations in
last 5 years (see \cite{AZ}, \cite{vD}, \cite{MD1}, \cite{Ner3}--%
\cite{Ner5}, \cite{NO}, \cite{OO}--\cite{OZ},
\cite{UU}, \cite{Zha}), and we will formulate  reasons
which are the most closed to the author.
 In first place, there are many explicit
 analytical formulas
related to kernel representations (I hope that this paper
also confirms this statement).
Secondly, spectra of kernel representations are very rich%
\footnote{ The most interesting spectral
 problems of noncommutative harmonic analysis
which were intensively investigated in last 20 years are

  --- $L^2$ on pseudoriemannian symmetric spaces

  --- Howe dual pairs
(and the problem of decomposition of
 $L^2$ on Stiefel manifolds which are in some sense equivalent to
Howe dual pairs)

 Each representation, which  occurs
in spectra of Howe dual pairs, occurs in spectra of some kernel
representation. Converse statement is false. The
{\tt a priori}
 explanation
of this phenomenon is contained in \cite{Ner3}.

 I think that  spectra
of kernel representations and spectra of $L^2$
on pseudoriemannian symmetric spaces
essentially differs. {\tt A priori} embedding
of spectra of $L^2(\U(p,q,\K)/\U(r,\K)\times\U(p-r,q,\K))$
for $\K=\R,\C,\Bbb H$
to spectra of kernel representations
is discussed in \cite{NO}, \cite{Ner3}.}.
Thirdly, the
kernel representations  have some interaction with
function theory%
\footnote{
For instance, in the work \cite{NO} we use functional theoretical
arguments for construction of singular unitary
representations of groups $\U(p,q,\K)$}, see \cite{NO}, \cite{Ner*},
\cite{Ner5}. Forth, the kernel representations
also are closely related to Olshanskii constructions
of representations of infinite dimensional groups
$\U(p,\infty)$, $\O(p,\infty)$, $\Sp(p,\infty)$
(see \cite{NO}).

\mal

{\bff 0.2. Bibliographical comments.}
Let $G/K$ be itself an hermitian symmetric space
 (i.e. $G=\U(p,q)$, $\Sp(2n,\R)$,
$\SO^*(2n)$, $\SO(n,2)$). Then its hermitization
$\widetilde G/\widetilde K$ is $G/K\times G/K$.
 A kernel representation of $G$
in this case is a tensor product of
 a highest weight representation $\widetilde\rho_\alpha$
of $G$
 and a lowest weight representation $\widetilde\rho_\beta^*$ of $G$.
 In  short paper \cite{Ber2} published in 1978
Berezin announced nice Plancherel formula
 for the {\it sufficiently large
parameter} $\alpha=\beta$ {\it of highest weight}
(see below Subsections 1.10-1.11).
In this case%
\footnote{Tensor products for ${\rm SL}_2(\R)=\U(1,1)$ were
 earlier investigated
by Pukanszky \cite{Puk} and by Vershik, Gelfand and
Graev \cite{VGG}, see also \cite{Mol}}
the kernel representation is
equivalent to the representation of $G$ in $L^2(G/K)$.
Berezin died soon after this and the proof never was published%
\footnote{ 20 years later I heard some reminiscences about this
proof but I can not reconstruct proof itself.
It  essentially differs from
 Unterberger--Upmeier \cite{UU}
proof and my proof \cite{Ner4}.
}.
Berezin work didn't attract serious interest in this time
 (see only papers \cite{Gut},\cite{Rep} on related subjects).

Secondly, the kernel representations
 appeared in G.I.Olshanskii and my work
which was partially announced in \cite{Ner}, \cite{Ols1}
 and partially published in \cite{NO}.
This work cocerned in vector valued kernel representations
of the groups $G=\O(p,q)$, $\U(p,q)$, $\Sp(p,q)$
for {\it small values} of the highest weight.
The main topic of our work
was investigation of discrete part
 of spectra of the kernel representations
and construction of ``exotic''  unitary representations
of $G$ by simple functional theoretical tools.

In the middle of 90-s interest to kernel representations
increases (we list some publications:
\cite{UU}, \cite{AZ}, \cite{Ner3}, \cite{OO},
\cite{OZ}, \cite{vD}, \cite{MD1},\cite{Zha}).
In 1994 Upmeier and Unterberger \cite{UU}
 published proof of Berezin formula
(see also \cite{AZ})\footnote{Their result also covers
groups $E_6$, $E_7$.}.
 Van Dijk and Hille \cite{vD} obtained the complete Plancherel formula for
rank 1 groups.
 Olaffson and Orsted  \cite{OO} proved that for a large highest weight
a scalar valued kernel representation of $G$ is equivalent to
the representation of $G$ in $L^2(G/K)$.

In  paper \cite{Ner4} there was defined
 $\B$-function for arbitrary classical
noncompact Riemannian symmetric space. For the symmetric
cones $\GL(n,\R)/\O(n)$, $\GL(n,\C)/\U(n)$,
 $\GL(n,\H)/\Sp(n)$ these $\B$-integrals coincides
with Gindikin $\B$-function constructed in \cite{Gin1}(1964)
(see also exposition in \cite{FK}). The construction
of \cite{Ner4} for  special cases of parameters
 gives some integrals
of Siegel \cite{Sie}, Hua Loo Keng \cite{Hua},
Unterberger--Upmeier \cite{UU}, Arazy--Zhang \cite{AZ}.
The Plancherel formula for
scalar-valued kernel
representations of all classical groups for large values of parameter
is easily reduced to these $\B$-integrals.
In this case a kernel representation is
equivalent to the representation of $G$ in $L^2(G/K)$
and the spectrum of the kernel representation
is supported by the principal nondegenerate unitary series.

The case of small values
of highest weight was discussed in paper \cite{Ner5}. In this case
the spectrum of kernel representation is quite intricate and  work
\cite{Ner5} contains
natural decomposition of a kernel representation on subrepresentations
having relatively simple spectra.

The purpose of the present parer is to obtain the complete
Plancherel formula for the kernel representations in the scalar-valued case.

\mal

{\bff 0.3. Contents.} Main part (Sections 1--5)
of the paper deals with the series
$G=\O(p,q)$\footnote{This case is the most complicated
and all difficulties  existing for other series exist
also for $\O(p,q)$. For all other series our proof is more simple.
}.

Section 1 of the paper contains preliminaries. We discuss
the definition of the  kernel representations
and simple a priori properties of the Plancherel formula. We also
formulate some necessary properties of spherical functions.

Basic results are formulated in Section 2. For large values of
$\alpha$ (where $\alpha$ is the
parameter of a highest weight) the Plancherel measure $\nu_\alpha$
has the form
$$E(\alpha)\prod_{k=1}^p |\Gamma(\frac12(\alpha+(p+q)/2+s_k)|^2 \r(s)ds$$
where $\r(s)$ is the Gindikin-Karpelevich density (see(1.42)--(1.43)),
\begin{equation}
s_1,\,\dots,\, s_p\in i\R
\end{equation}
(this notation means that $\Re s_j=0$)
and $E(\alpha)$ is a meromorphic factor.

Assume  $q-p$ be sufficiently large.
Let us move the parameter
 $\alpha$ from $+\infty$ to 0. After passing
across the point $\alpha=\frac12(p+q)-1$
there appears an additional piece of the support
 of the Plancherel measure. This piece is defined by the
conditions
\begin{equation}
s_1=\alpha- {\t\frac12}({p+q})+1;\quad s_2,\dots,  s_p\in i\R
\end{equation}
After passing across the point $\alpha=\frac12(p+q)-2$
  the third piece of the support of the Plancherel measure appears:
\begin{equation}
s_1=\alpha-\half(p+q)+1,\,s_2=\alpha-\half(p+q)+2 ;\quad
        s_3,\dots, s_p\in i\R
\end{equation}
After passing across the point $\alpha=\frac12(p+q)-3$
we  obtain two additional components of the support
\begin{multline*}
s_1=\alpha-\half(p+q)+1,\, s_2=
\alpha-\half(p+q)+2,\,s_3=
\alpha-\half(p+q)+3;\quad s_4,\dots,  s_p\in i\R
\end{multline*}
and
\begin{equation}
s_1=\alpha-\half(p+q)+3;\quad s_2,\dots,  s_p\in i\R
\end{equation}
etc. After passing across the point $\alpha=(q-p)/2$
 we  obtain the first one-point piece
$$s_1=\alpha-\half(p+q)+1;\quad s_2=
\alpha-\half(p+q)+2,\dots, s_p=
\alpha-\half(p+q)+p
$$
 This means that our
representation has a subrepresentation entering discretely.

At the point $\alpha=p-1$ the component (0.1) of the support
 disappears. At the point $\alpha=p-2$
components (0.2), (0.4) also disappear, etc.

Theorems 2.2--2.4  contain the complete description of this process and
give the Plancherel density on each component of the support.
Interpretation of these pieces is given in \cite{Ner5},
in the present paper this is not discussed.

The nature of spectra of kernel representations
is explained in \cite{Ner5}%
\footnote{The discrete part of spectra
of kernel representations consists
of singular unitary representations having quite interesting
properties. For instance, these infinite dimensional
(non highest weight) representations have Gelfand--Tsetlin bases, see
\cite{Ols1}, \cite{Molev}; the problems of decomposition
of restrictions and tensor products for these representations
also seems rich, see \cite{Ner}, \cite{NO}).
In \cite{Ols1}, \cite{NO} it was shown, that
these representations  admit inductive limits as $q\to\infty$.
 Certain representations of
this type appear in spectra discussed in \cite{kru}, \cite{Kob}}.

For an integer negative $\alpha$ our
 construction gives the Plancherel formula
for some finite dimensional representation of $\O(p,q)$
(see Subsection 2.6 of the paper).

Section 3 is based on \cite{Ner4} and  contains
 evaluation of the $\B$-integral (see formula (3.2)-(3.4)).
 For instance, in the  case
$p=q$ our $\B$-integral is given by

 \begin{align}&
\int_{R+R^t>0}
\prod_{j=1}^p
\frac{ \det[(R+R^t)/2]_j^{\lambda_j-\lambda_{j+1}}  }
  {\det [1+R]_j^{\sigma_j-\sigma_{j+1}} }
  \cdot
\det(R+R^t)^{-p} dR
=\\&
=\const  \cdot
\frac{\Gamma(\lambda_k-(p+k)/2+1)\Gamma(\sigma_k-\lambda_k-(p-k)/2)}
     {\Gamma(\sigma_k-p+k)}    \nonumber
\end{align}
where the integration is given over the
 space of dissipative $p\times p$ real matrices
$R$ and the symbol $[A]_j$ denotes the left upper $j\times j$
block of a matrix $A$.

$\B$-Integral allows to obtain the
 Plancherel formula for $\alpha>\frac12(p+q)-1$.
In Section 4 we construct the analytic continuation of the Plancherel formula
to  arbitrary $\alpha$.

In Section 5 we prove positive
 definiteness of spherical functions
which appears in the right side of the Plancherel formula.

Section 6 contains a discussion of other series of classical groups.
The $\B$-integrals for other series or real classical groups
 are evaluated in
\cite{Ner4}, and a generalization of the consideration
of Subsections 1, 4, 5 to other series is quite trivial. Hence, we give
only short remarks and also give the Plancherel formula
in the form which slightly differs from Theorem 2.2.
Author intentionally considers the series $\O(p,q)$
(and not so-called 'general case')
to do the exposition
 more or less self-closed.
I try
to avoid formal logical dependence on recent papers and
also  minimize using machinery of representation
theory of semisimple groups as far as it is possible%
\footnote{we needs in some basic properties of spherical
functions, all necessary information is contained
in Helgason book \cite{Hel}, chapter 4};
I also try to avoid notations demanding  long explanations.

\mal

{\bf Acknowledgements.}
I am very grateful to G.I.Olshanskii,
V.F.Molchanov, and B.Orsted for numerous
discussions of the subject. I
thanks H.Schlichtkrull, G. van Dijk,
and A.Dvorsky for discussions,
 comments and references.

\bol

{\bfff  1. Preliminaries.}

\nopagebreak

\bol

\addtocounter{sec}{1}
\setcounter{equation}{0}

{\bfff  A. Positive definite kernels.}

\nopagebreak

\bol

The subject of the paper is the analysis in a family
of hilbert spaces defined by positive definite kernels.
 The notion of positive definite kernel and
associated machinery are quite old
(see  \cite{Sch}, \cite{Berg}, \cite{Kre})
but not widely known.
In this Section we briefly discuss elementary properties of
the positive definite kernels and associated hilbert spaces.

{\bff 1.\punkt. Positive definite kernels}.
Let $H$ be a hilbert space with a scalar product
$<\cdot,\cdot>$, let $X$ be a subset in $H$.
Consider the function $L(x,y)$ on $X\times X$
defined by
$$L(x,y)=<x,y>$$
Obviously for all $x_1,\dots,x_n\in X$ we have
\begin{equation}
\det\begin{pmatrix}
 L(x_1,x_1)&\cdots&L(x_1,x_n)\\
 \cdots&\cdots&\cdots\\
 L(x_n,x_1)&\cdots&L(x_n,x_n)
\end{pmatrix}\ge 0
\end{equation}

 Let $X$ be an abstract set.
A   function $L(x,y)$ on $X\times X$ is called
a {\it positive definite  kernel} if it satisfies the conditions

1. $L(x,y)=\ov{L(y,x)}$

2. For any $x_1, x_2,\dots, x_n\in X$ inequality (1.1) holds

Let $L(x,y)$ be a positive definite kernel on $X$. Then
where exists a  hilbert space $H=H[K]$ and a system of vectors
$v_x\in H$ enumerated by points $x\in X$ such that

1. $<v_x,v_y>_H=L(x,y)$

2. the linear span of the vectors $v_x$ is dense in $H$.

The family $v_x$ is called a  {\it supercomplete basis}\footnote{
Other terms for $v_x$ are {\it overfilled basis} or {\it system
of coherent states}.}.

This construction is natural in the following sense.
Let $H'$ be another hilbert space
and let $v_x'$ be another system of vectors satisfying
the same conditions. Then there exists the unique unitary operator
$U:H\to H'$ such that $Uv_x=v_{x}'$ for all $x\in X$.

If $X$ is a separable metric space and the kernel $L(x,y)$ is continuous,
then the hilbert space $H[L]$ is separable.

In Subsections 1.2-1.3 and 1.4
 we  discuss two ways of ``materialization'' of the space $H[L]$.

\mal

{\bff 1.\punkt. Scalar product
in the space of complex-valued measures.}
Assume  $X$ be a separable complete metric space.
Let $\mu$ be a complex-valued measure
(charge) on $X$ with a compact support.
Consider a  vector
$$v(\mu)=\int_X v_x  d\mu(x)\in H[L]$$
Thus, we obtain a way to represent  elements of $H[L]$.
Obviously
$$<v(\mu), v(\nu)>_{H[L]}=\int_{X\times X} L(x,y)\,d\mu(x)\,d\ov{\nu(y)}$$

Let us say the same construction more formally.
 Consider the linear space ${\cal M}(X)$
of all compactly supported complex-valued
 measures  on $X$.
Consider the scalar product in ${\cal M}(X)$ defined by the formula
$$<\mu, \nu>=\int_{X\times X} L(x,y)\,d\mu(x)\,\ov{d\nu(y)}$$
We obtain a structure of a prehilbert space in ${\cal M}(X)$
and the space $H[L]$ is the  hilbert space associated with the prehilbert
space   ${\cal M}(X)$
(elements of the supercomplete basis corresponds to
measures supported at points).

\mal

{\bff 1.\punkt. Scalar products in  spaces of distributions.}
Assume  $X$ be a smooth manifold and the kernel $L(x,y)$ be smooth.
Denote by ${\cal D}$ the space of compactly supported distributions
on $X$.
 Consider the scalar product in ${\cal D}$
given by the formula
\begin{equation}
<\chi,\phi>=\{K(x,y),\chi(x)\otimes \ov{\phi(y)}\}
\end{equation}
where brackets
 $\{\cdot,\cdot\}$ denote the pairing of  smooth functions and
distributions. Consider the hilbert space $H$
 associated with the prehilbert space ${\cal D}$.
Denote $\delta$-distribution supported at
a point $x$ by $\delta_x$. Obviously
$$<\delta_x,\delta_y>=L(x,y)$$
Hence, we can identify $H$ with $H[L]$ and the vectors
$\delta_x$ with elements of the supercomplete basis $v_x$.

     \mal

{\scc Remark 1.} The space of distributions equipped with  scalar product
(1.2) is not complete. This means that some vectors of
$H[L]$  can not be represented by distributions.

\mal

{\scc Remark 2.} Scalar product (1.2) in ${\cal D}$
  can be degenerated. This means
that a vector $h\in H$ can be represented by a distribution in various ways.

\mal

{\bff 1.\punkt. The embedding of $H[L]$ to the space of functions on $X$.}
For arbitrary $h\in H[L]$ we consider the function
$$f_h(x):=<h,v_x>_{H[L]}$$
on the space $X$. Obviously, the map $h\mapsto f_h$
is an embedding of $H[L]$ to the space of functions on $X$.
We denote the image of the embedding by
$H^\circ[L]$. By construction, the space $H^\circ[L]$
has structure of a hilbert space.

\mal

{\scc Lemma 1.\fact.} {\it Assume  $X$ be a separable metric space and
the kernel $L(x,y)$ be continuous. Let $h_j$ converges to $h$. Then
$f_{h_j}$ converges to $f_h$ uniformly on compacts.}

\mal

{\scc Proof.} Let $Y\subset X$ be a compact set.
Let $x\in Y$. Then
\begin{multline*}
 |f_{h_j}(x)-f_h(x)|=|<(h_j-h, v_x>_{H[L]}|\le
\\ \le
 \|h_j-h\|\cdot\|v_x\|=
 \|h_j-h\|\cdot \sqrt{L(x,x)}\le
 \|h_j-h\|\cdot \sqrt{\max\limits_{x\in Y} L(x,x)}
\qquad\boxtimes
\end{multline*}

 Obviously, the function $\phi_a(x)\in H^\circ[L]$ associated
with the vector
$v_a\in H[L]$ is given by the formula
\begin{equation}
\phi_a(x)=L(x,a)
\end{equation}

{\scc Lemma 1.\fact.} {\it  Let $X$ be a locally compact metric space
and the kernel $L(x,y)$ be continuous.
 Then  functions $f_h\in H^\circ[L]$ are continuous.}

\mal

{\scc Proof.} The linear span of the functions $\phi_a$ is dense
in $H^\circ[L]$. Then we apply
Lemma 1.1.        \hfill  $\boxtimes$

\mal

{\scc Lemma 1.\fact.} ({\sc Reproducing property}) {\it
For any $f\in H^\circ[L]$, $x\in X$ the following identity holds}
\begin{equation}
f(x)=<f,\phi_x>_{ H^\circ[L]}
\end{equation}

{\scc Proof.} Let $f=f_h$. Then
$$<f_h,\phi_x>_{ H^\circ[L]}= <h, v_x>_{H[L]}=f_h(x)
\qquad\qquad\qquad\boxtimes$$

 {\scc Remark.}
Equation (1.4) gives a nonexplicit description of the scalar product
in $  H^\circ[L]$ and this description
 is sufficient for many purposes.
Another way of description  is the following
identity.
Let $e_n(x)\in H^\circ[L]$ be an orthonormal basis. Then
$$
L(x,y)=\sum\ov{e_n(x)}e_n(y)$$
(proof: consider $<\phi_x, e_n>$)

\mal

{\scc Lemma 1.\fact.} {\it  Let $\Omega$ be an open domain in $\C^n$ and
let $L(z,u)$ be a positive definite kernel on $\Omega$.
Let $L(z,u)$ be holomorphic in the variable $u$
 and anti-holomorphic in  the variable
$z$. Then all elements of
 the space $H^\circ[L]$ are holomorphic functions
on $\Omega$.}

\mal

{\scc Proof.} It is a corollary of Lemma 1.1. \hfill $\boxtimes$

\mal

{\bff 1.\punkt. Operations with positive definite kernels.}

{\scc Lemma 1.\fact.}
a) {\it Let $L_1(x,y)$, $L_2(x,y)$ be positive definite kernels on $X$.
Then $L_1(x,y)+L_2(x,y)$ is a positive definite kernel}

b) {\it Let $L_1(x,y)$, $L_2(x,y)$ be positive definite kernels on $X$.
Then the kernel $L_1(x,y)L_2(x,y)$ is  positive definite}

c) {\it Let $L_j(x,y)$ be positive definite kernels
 and $L_j(x,y)$ converges
to $L(x,y)$ point-wise. Then $L(x,y)$ is a positive definite
kernel.}

d) {\it Let $K_m(x,y)$ be a family of positive definite kernels
 enumerated by points of some measure space $M$
with positive measure $\mu$. Assume that the integral
$$K^*(x,y)=\int_M K_m(x,y)d\mu(m)$$
converges for all $x,y\in X$. Then $K^*(x,y)$
is positive definite.}

e) {\it Let $L(x,y)$ be a positive definite kernel and let
$\lambda(x)$ be a function on $X$. Then the kernel
$M(x,y)=\lambda(x)\ov{\lambda(y)}L(x,y)$ is positive definite.}

{\scc Proof.} a) Let $v_x$ (resp. $w_x$) be the supercomplete basis
in $H[L_1]$ (respectively $H[L_2]$). We consider the system of vectors
$v_x\oplus w_x\in H[L_1]\oplus H[L_2]$. Then
$$L_1(x,y)+L_2(x,y)=<v_x\oplus w_x,v_y\oplus w_y >$$

b) Proof is similar,
$H[L_1L_2]\subset H[L_1]\otimes H[L_2]$

d) This is consequence of
a) and c).

e) Indeed, $H[M]=H[L]$ and the supercomplete basis in $H[M]$ consists
of vectors $\gamma(x)v_x$ where $v_x$ is the supercomplete basis in
$H[L]$.\hfill$\boxtimes$

\mal

{\bff 1.\punkt. Positive definite kernels on homogeneous spaces.}
Let $\Gamma$ be a group acting on $X$ and let a
positive definite kernel $L(x,y)$
be $\Gamma$-invariant
$$L(gx,gy)=L(x,y) \qquad\mbox{for all}\qquad g\in \Gamma, x,y\in X$$
Obviously,
 for each $g\in\Gamma$ where exists
the unique unitary operator
$U(g):H[L]\to H[L]$ such that
$$Uv_x=v_{gx}\qquad \mbox{for all} \qquad x\in X$$
Then
$$U(g_1g_2)=U(g_1)U(g_2)$$
Hence, $U(g)$ is an unitary representation of $\Gamma$.

Let $X$ be a $\Gamma$-homogeneous space, $X=\Gamma/K$,
 let $x_0$ be a $H$-fixed
point.

 Let $L(x,y)$ be a $\Gamma$-invariant function. Then $L$ is
completely defined by the function
$$l(y):=L(x_0,y)$$
 Indeed, let $u,z\in X$.
Then $u=gx_0$ for some element $g$ in $\Gamma$ and
$$L(u,z)=L(gx_0,z)=L(x_0,g^{-1}z)=l(g^{-1}z)$$
Moreover, for any $\gamma\in K$ we have
$$l(y)=L(x_0,y)=L(\gamma x_0, \gamma y)=L(x_0,\gamma y)=l(\gamma y)$$
We see that the function $l_0$ is a $K$-invariant function on
$\Gamma/K$.

We also can consider a $K$-invariant function $l(y)$ as
function on double cosets $K\setminus \Gamma/K$.

We see that {\it there is the canonical correspondence between 3
following sets:

-- $\Gamma$-invariant functions on $\Gamma/K\times \Gamma/K$

-- $K$-invariant functions on $\Gamma/K$

-- functions on $K\setminus \Gamma/K$.}

We say that a $K$-invariant function on $\Gamma/K$ or
a function  on $K \setminus \Gamma/K$ is
{\it positive definite} if the associated kernel on
$\Gamma/K\times \Gamma/K$ is positive definite.

\mal

{\bff 1.\punkt. On $K$-invariant vectors in representations
of $\Gamma$}. Let $\Gamma$, $K$, $x_0$ be the same as above.
 Let $\rho$ be an unitary representation  of $\Gamma$
in a hilbert space $H$. Assume that there exists a
$K$-invariant vector $v\in H$ and assume  $v$ be
a {\it cyclic vector}%
\footnote{This means that the linear span of vectors
$\rho(g)v$ is dense in $H$.}.

Consider the map $\Gamma/K\to H$ given by the formula
$$gx_0\mapsto \rho(g)v$$
(the image of the map is the $\Gamma$-orbit of the vector $v$).
Then the    function
$$L(g_1x_0,g_2x_0):=<\rho(g_1)v, \rho(g_2)v>_H$$
is a $\Gamma$-invariant positive definite kernel on $\Gamma/K$.

 Hence, we can identify  the hilbert space $H$ with the space
$H[L]$; the supercomplete basis in  $H$ consists of vectors
$\rho(g)v$.

 The function $l(y)$ in our case is the matrix element
$<\rho(g)v,v>$  and our construction
(Segal--Gelfand--Naimark construction) reconstructs the representation
$\rho$ by its matrix element.

\bol

{\bfff B. Kernel representations.}

\nopagebreak

\bol

Assume $p\le q$

{\bff 1.\punkt. Pseudoorthogonal group $\O(p,q)$.} Consider the linear space
$\C^p\oplus\C^q$ equipped with the indefinite hermitian form
\begin{equation}
J((x,y),(u,v))=\sum_{j=1}^p x_j\ov{u}_j- \sum_{j=1}^q y_k\ov v_k;\qquad
(x,y),(u,v)\in \C^p\oplus\C^q
\end{equation}
The {\it pseudounitary group} $\U(p,q)$ is the group of all linear operators
$g=\begin{pmatrix}\alpha&\beta\\ \gamma&\delta\end{pmatrix}$ in
$ \C^p\oplus \C^q$ preserving the form $J(\cdot,\cdot)$.
In other words,
 a matrix $g\in \U(p,q)$ satisfies the condition
\begin{equation}
\begin{pmatrix}\alpha&\beta\\ \gamma&\delta\end{pmatrix}
\begin{pmatrix}1&0\\0&-1\end{pmatrix}
\begin{pmatrix}\alpha&\beta\\ \gamma&\delta\end{pmatrix}^*=
\begin{pmatrix}1&0\\0&-1\end{pmatrix}
\end{equation}
The {\it pseudoorthogonal group} $\O(p,q)$ is the subgroup of
$\U(p,q)$ consisting of real matrices. Below in Sections 1-5
 by the symbol
$G$ we denote the group
$$G=\O(p,q)$$
By $K$ we denote the subgroup $\O(p)\times\O(q)\subset G$
consisting of matrices
having the form $\begin{pmatrix}\alpha&0\\0&\delta\end{pmatrix}$.
It is a maximal compact subgroup in $G$.

{\bff 1.\punkt. Matrix balls.} By $\B_{p,q}(\C)$ we denote space of
all complex $p\times q$ matrices $z$ having norm $<1$
(where {\it a norm} is the norm of the operator $v\mapsto vz$
from the euclidean space $\C^p$ to the euclidean space $\C^q$;
remind that $\|z\|^2$ is the maximal eigenvalue of $z^*z$).

By $\B_{p,q}(\R)$ we denote the space of real $p\times q$ matrices
with norm $<1$.

The group $\U(p,q)$ acts on the matrix ball $\B_{p,q}(\C)$ by
{\it fractional linear}
transformations
\begin{equation}
z\mapsto z^{[g]}:=(\alpha+z\gamma)^{-1}(\beta+z\delta)
\end{equation}
This action is transitive and
the stabilizer of the point $z=0$ is the subgroup $\U(p)\times\U(q)$.
Hence, $\B_{p,q}(\C)$ is the symmetric space
$$\B_{p,q}(\C)=\U(p,q)/\U(p)\times\U(q)$$
In the same way, $\B_{p,q}(\R)$ is the symmetric space
$$\B_{p,q}(\R)=G/K=\O(p,q)/\O(p)\times\O(q)$$
Arbitrary symmetric space admits unique up to factor invariant measure.
For the space $\B_{p,q}(\R)$ the $\O(p,q)$-invariant measure
is given by the formula
\begin{equation}
d\lambda(z)=\det(1-z^*z)^{-(p+q)/2}d\mu(z)
\end{equation}
where $d\mu(z)$ is the Lebesgue measure on $\B_{p,q}(\R)$.

\mal

{\bff 1.\punkt. Berezin kernels.}
 {\scc Theorem 1.\fact\footnote%
{see Berezin \cite{Ber1} (1975), see also Gindikin \cite{Gin1},
Rossi, Vergne \cite{RV}, Wallach \cite{Wal},
see also a recent exposition in \cite{FK}.} }.
{\it The kernel
$$L_\alpha(z,u)=\det(1-z^*u)^{-\alpha}$$
on the matrix ball $\B_{p,q}(\C)$ is positive definite if and only if}
\begin{equation}
\alpha=0,1,2,\dots,p-1\quad\mbox{or}\quad \alpha>p-1
\end{equation}

Thus, for $\alpha$ satisfying the Berezin condition
(1.9), we obtain the hilbert spaces $H_\alpha:=H[L_\alpha]$
and $H_\alpha^\circ:=
H^\circ[L_\alpha]$.
The function $L_\alpha(z,u)$ is anti-holomorphic in $z$
and hence by Lemma 1.4 the space $H_\alpha^\circ$
consists of holomorphic functions on  the matrix
ball $\B_{p,q}(\C)$.

\mal

{\scc Remark.} For $\alpha>p+q-1$
the scalar product in $H_\alpha^\circ$ can be represented in the form
$$<f,g>_\alpha=
   C(\alpha)\int_{\B_{p,q}}f(z)\ov{g(z)}\det(1-z^*z)^{\alpha-p-q}d\mu(z)$$
where $d\mu(z)$ is the  Lebesgue measure on $\B_{p,q}(\C)$
and $C(\alpha)$ is the meromorphic factor defined by the condition
$<1,1>_\alpha=1$. In particular
 $H_{p+q}^\circ$ is the Bergman space.
For $\alpha=q$ we obtain the Hardy space $H^2$.
The scalar product in this case
is given by the formula
$$
<f(z),g(z)>_q=\int_{z^*z=1}f(z)\ov{g(z)}d\nu(z)
$$
where $d\nu(z)$ is the unique $\U(p)\times\U(q)$-invariant
measure on the set\footnote{Schtiefel manifold}
$z z^*=1$. For other values of parameters there exist
integral formulas including partial derivatives but they are not simple
(see \cite{AU}).

\mal

{\scc Remark.}
For $\alpha>p-1$ the space $H_\alpha$ contains all polynomials
on $\B_{p,q}$. For $\alpha=0,1,\dots,p-1$ all functions $f\in H_\alpha$
satisfy  some system of partial differential equations.
For $\alpha=0$ our space contains only constants.

\mal

{\scc Proposition 1.\fact.} (\cite{Ber1}, \cite{RV})
 a){\it For any $g=
\begin{pmatrix}a&b\\ c&d\end{pmatrix}\in \U(p,q)$
the operator
\begin{equation}
\widetilde T_\alpha(g)
(z)f(z)=f((a+zc)^{-1}(b+zd))\det(a+zc)^{-\alpha}
\end{equation}
is unitary in $H^\circ_\alpha$}

{\scc Remark.} If $\alpha$ is not integer, then
\begin{multline}
\det(a+zc)^{-\alpha}=\det a^{-\alpha}\det(1+zca^{-1})^{-\alpha}=\\
=|\det a|^{-\alpha}\cdot e^{-\alpha (i\arg \det a+ 2\pi ki)}
\det(1+zca^{-1})^{-\alpha}
\end{multline}
 is a multi-valued function. It is easy to show that
$\|ca^{-1}\|<1$. Hence, $(1+zca^{-1})^{-\alpha}$
is a well-defined single-valued function on the matrix ball $\B_{p,q}(\C)$.
Hence, expression (1.11) has countable family of holomorphic branches
 on $\B_{p,q}(\C)$ and   formula (1.10)
defines a countable family of well-defined operators
which differs by constant factors $e^{2\pi k\alpha i}$.

\mal

{\scc Proof.} Consider the supercomplete basis
 $\phi_x(z)=\det(1-x^*z)^{-\alpha}$ in $H_\alpha^\circ$.
 A calculation show that
$$\widetilde T_\alpha(g^{-1}) \phi_x(z)=
        \det(a+zc)^{\alpha}\phi_{x^{[g]}}(z)$$
The simple identity
 $$ \det(1-x^{[g]} (y^{[g]})^*)=
\det(1-xy^*)\det(a+xc)^{-1} \det(a+yc)^{-1}$$
implies
$$<\phi_x   , \phi_y>_{H^\circ_\alpha}=
<\widetilde T_\alpha(g^{-1})\phi_x   ,
\widetilde T_\alpha(g^{-1}) \phi_y>_{H^\circ_\alpha}
           \qquad\qquad\qquad\boxtimes
$$

Obviously
$$\widetilde T_\alpha(g_1) \widetilde T_\alpha( g_2)
=e^{2\pi m\alpha i}\widetilde T_\alpha(g_1 g_2),
     \qquad\mbox{where} \quad m\in {\Bbb Z}$$

If $\alpha$ is integer, then $\widetilde T_\alpha$ is a linear
representation of $\U(p,q)$.
If $\alpha$ is not integer, then $\widetilde T_\alpha$ is
a projective representation of $\U(p,q)$ or
a linear representation of the universal covering
group $\U(p,q)^\sim$ of the group
$\U(p,q)$.

\mal

{\bff 1.\punkt. Kernel representations of $\O(p,q)$.}
 {\it  Kernel representation}
$T_\alpha$ of the group $G=\O(p,q)$ is the restriction of
the representation $\widetilde T_\alpha$ to the subgroup $\O(p,q)$.
We also say that the function $f(z)=1$ is {\it the marked vector}
in $H^\circ_\alpha$. We denote this vector by $\Xi$.

\mal

{\scc Remark.} A kernel representation is a linear representation.
Indeed, we can wright $|\det a|^{-\alpha} (1+zca^{-1})^{-\alpha}$
instead of (1.11).

\mal

{\scc Lemma 1.\fact.} {\it The vector
$\Xi$ is $\O(p,q)$-cyclic.}
% \rm(\it i.e. linear span of $\O(p,q)$-orbit
%of $\Xi$ is dense in $H^\circ_\alpha$.\rm)}

{\scc Proof.} Let $Q\subset H_\alpha$ be a subspace containing
the $G$-orbit of $\Xi$. This orbit  consists of functions
(1.11) and hence the functions $\det(1+zca^{-1})^{-\alpha}$ are contained in $Q$.
But the point $ca^{-1}$ is the image of $0$ under
the fractional linear
 transformation (1.7). Since the action of $\O(p,q)$ is transitive
on $\B_{p,q}(\C)$, the subspace $Q$ contains
all functions $v_u=\det(1+z u^*)^{-\alpha}$ where
$u\in\B_{p,q}(\R)$. Furthermore,                      since
the  family $v_u$ depends on $u$ holomorphically, $v_u\in Q$
for all $u\in\B_{p,q}(\C)$. But $v_u$ is the supercomplete basis in
$H^\circ_\alpha$. Hence, $Q=H^\circ_\alpha$.  \hfill$\boxtimes$

\mal

{\scc Lemma 1.\fact.}
{\it Any $\O(p,q)$-invariant subspace in $H_\alpha^\circ$ contains an
 $\O(p,q)$-invariant vector.}

{\scc Proof.} Assume $H_\alpha^\circ=R\oplus Q$ where $R$, $Q$ are
invariant subspaces. Assume that $R$ hasn't an $\O(p)\times\O(q)$-%
invariant vector. Then the projection of $\Xi$ to $R$ is zero, and hence
$\Xi\in Q$. But $\Xi$ is cyclic. Thus, $Q=H^\circ_\alpha$.
            \hfill$\boxtimes$

\mal

{\bff 1.\punkt. Another description of the kernel-representations.}
Let $\alpha$ satisfies Berezin conditions
(1.9). By Lemma 1.5.e)
the kernel
\begin{equation}
M_\alpha(z,u)=\frac{\det(1-zz^*)^{\alpha/2}\det(1-uu^*)^{\alpha/2}}
                     {\det(1-zu^*)^\alpha}
\end{equation}
on $\B_{p,q}(\R)$ is positive definite.
 A simple calculation show that the
kernel $M_\alpha$ is $\O(p,q)$-invariant.
Hence, we obtain an unitary representation of the group $\O(p,q)$
in the hilbert space $H[M_\alpha]\simeq H^\circ[M_\alpha]$
(see Subsection 1.6).
The group $\O(p,q)$ acts in $H^\circ[M_\alpha]$ by substitutions
\begin{equation}
f(z)\mapsto f((a+zc)^{-1}(b+zd))
\end{equation}

The {\it marked vector} $\Xi$
in this model is the element of the supercomplete basis corresponding
to the point $0\in \B_{p,q}(\R)$.

 Let us define the canonical unitary
$\O(p,q)$-intertwining operator
 $$A:H^\circ[L_\alpha]\to H^\circ[M_\alpha]$$
 Let $f\in H^\circ[L_\alpha]$, let $z\in \B_{p,q}(\R)$. Then
$$Af(z)=f(z)\det(1-zz^*)^{\alpha/2}$$

This map transforms elements of the supercomplete basis
in $H^\circ[L_\alpha]$ to elements of the supercomplete basis
in $H^\circ[M_\alpha]$.

{\bff 1.\punkt. Limit as $\alpha\to\infty$.}
Let $\lambda$ be the $\O(p,q)$-invariant measure
on $\B_{p,q}(\R)$ (see (1.8)). Denote by $C_0$
 the space
of continuous functions on $\B_{p,q}(\R)$ with a compact support.
If $\phi\in C_0$, then $\phi(z)\lambda(z)$ is a complex valued
 measure on $\B_{p,q}$. Hence (see Subsection 1.2),
we obtain the scalar product in the space  $C_0$
given by
\begin{equation}
<\phi,\psi>=A_\alpha \int_{\B_{p,q}(\R)\times\B_{p,q}(\R)}
 M_\alpha(z,u) \phi(z)\ov{\phi(u)}
\,d\lambda(z) d\lambda(u)
\end{equation}
Let us define the normalization constant $A_\alpha$ by the
condition
$$A_\alpha=\Bigl(\int_{\B_{p,q}(\R)}(1-zz^*)^\alpha dz\Bigr)^{-1}$$
(it is a Hua Loo Keng integral, see (3.5)).
Obviously $M_\alpha(z,z)=1$
and $M_\alpha(z,u)<1$ if $z\ne u$. It is easy to see that
 the sequence
$A_\alpha M_\alpha(z,u)$ approximates the distribution
$\delta(z-u)$. Thus, the limit of scalar products
(1.14) as $\alpha\to\infty$ is
$$ <\phi,\psi>= \int_{\B_{p,q}(\R)}
  \phi(z)\ov{\psi(z)}
\,d\lambda(z)$$
In this sense {\it the limit of kernel representations as
$\alpha\to\infty$ is the space $L^2(G/K)$}.

We emphasis that the action of $\O(p,q)$
in $L^2(\O(p,q)/\O(p)\times\O(q))$ and in all
spaces $H^\circ[M_\alpha]$ is given by
the same formula (1.13) and only scalar product
in the space of functions varies. We will see
that the spectrum of the representation $T_\alpha$
 and the structure of Plancherel
formula
essentially depends on $\alpha$.

\mal

{\bff 1.\punkt. Preliminary remarks on
the Plancherel formula.} Our purpose
is to obtain a decomposition of the kernel representation
$T_\alpha$ on irreducible representations.

 An irreducible representation of $G=\O(p,q)$ is called {\it spherical} if it
contains a $K$-fixed vector.
 This vector is called
 {\it spherical vector}. Remind that the space
of $K$-fixed vectors for $G$ has dimension 0 or 1
(Gelfand theorem, see for instance \cite{Hel}, Theorem 4.3.1
and Lemma 4.3.6). Denote the set of all
unitary spherical representations
of $\O(p,q)$ by $\widehat G_{sph}$

\mal

{\sc Remark.} The explicit description of this set
is not known.
Parametrization of all (generally speaking nonunitary)
 spherical representations
of $\O(p,q)$ is simple and it is given below in Subsection 1.17).

\mal

 By $H_\rho$ we denote the space of
a spherical representation $\rho$,
by $\xi(\rho)$ we denote the spherical vector
in $H_\rho$ whose length
is 1.

\mal

{\scc Lemma 1.\fact.} {\it Decomposition
of the kernel-representation $T_\alpha$  has the form
\begin{equation}
T_\alpha(g)=\int_{\rho\in\widehat G_{sph}} \rho(g)d\nu_\alpha(g)
\end{equation}
where $\nu_\alpha$ is a Borel measure on $\widehat G_{sph}$.}

\mal

{\scc Remark.} For the definition of direct integrals of representations
and the abstract Plancherel formula see, for instance, \cite{Kir}, 8.4.

\mal

{\scc Proof.} By Lemma 1.9, the decomposition contains only
spherical representations. Hence, by the abstract Plancherel theorem
the representation $T_\alpha(g)$ has the form
$$T_\alpha(g)= \bigoplus_{j=1}^{\kappa}R_j$$
where
$$R_j=\int_{\rho\in\widehat G_{sph}} \rho(g)d\nu_\alpha^{j}(g)$$
and  the measure $\nu^{j+1}_\alpha$ is absolutely continuous
with respect to $\nu^j_\alpha$ for all $j$. The number
 $\kappa$  can be 1,2,\dots, $\infty$.
We must  prove that $\kappa=1$.

All $K$-fixed vectors in $R_j$
are functions having the form
$\phi_j(\rho)\xi(\rho)$ where $\phi_j(\rho)$ is a $\nu^1_\alpha$-measurable
function on $  \widehat G_{sph}$.

Consider the projection $\Xi^{(1,2)}$ of the marked vector $\Xi$ to
$R_1\oplus R_2$. Since the vector $\Xi $ is cyclic in whole space,
 its projection must  be cyclic in $R_1\oplus R_2$.
 The vector $\Xi^{(1,2)}$ has the form
$$(\phi_1(\rho)\xi(\rho),\phi_2(\rho)\xi(\rho))\in R_1\oplus R_2 $$
Obviously the cyclic span of $\Xi^{(1,2)}$ in $R_1\oplus R_2$ contains
only vectors $$(q_1(\rho)\xi(\rho),q_2(\rho)\xi(\rho))$$
 satisfying the condition
$$\phi_2(\rho) q_1(\rho)=\phi_1(\rho) q_2(\rho)$$
If $\nu^2_\alpha\ne0$ we obtain a contradiction, since
the cyclic span of $\Xi^{(1,2)}$
 is a proper subspace in $R_1\oplus R_2$.   \hfill$\boxtimes$

\mal

{\bff 1.\punkt. Normalization of the Plancherel measure.}
The measure $\nu_\alpha$ in (1.15) is defined up to equivalence
of measures%
\footnote{measures $\mu$, $\nu$ are equivalent if there exists
a function $\chi$ such that $\chi\ne0$ almost everywhere
(in sense of $\nu$)
 and
$\mu=\chi\nu$.}.

The image of the marked vector $\Xi$ in the direct integral
 (1.15) is some function
$\phi(\rho)\xi(\rho)$ where $\xi(\rho)$ is an unit $K$-fixed
vector in $H_\rho$. It is convenient to assume
\begin{equation}
\phi(\rho)=1
\end{equation}
This assumption uniquely defines the measure $\nu_\alpha$.

\mal

{\sc Remark.} Assumption (1.16) is not restrictive. Indeed, let
us assume
that the image of $\Xi$ in (1.15)
is a function $\gamma(\rho)\xi(\rho)$. Then the Plancherel
measure is completely defined by this assumption
and it equals to
$\frac1{\sqrt {|\gamma|}}\nu_\alpha$.

After normalization (1.16) we obtain the following equality
of matrix elements
\begin{equation}<T_\alpha(g)\Xi,\Xi>_{H_\alpha}=<\left[\int_{
\widehat G_{sph}}\rho(g)d\nu_\alpha(\rho)\right]\cdot 1, 1>
\end{equation}
or
\begin{equation}
<T_\alpha(g)\Xi,\Xi>_{H_\alpha}=
\int_{\widehat G_{sph}}
 <\rho(g)\xi(\rho),\xi(\rho)>_{H_\rho}d\nu_\alpha(\rho)
\end{equation}
Conversely, assume that we know a measure $\nu_\alpha$
 on $\widehat G_{sph}$
satisfying  condition (1.18). Then it satisfies  condition (1.17).
Hence, the representations
 $T_\alpha$ and $\int_{\widehat G_{sph}}\rho(g)d\nu_\alpha$
have
the same matrix elements, and therefore they are
canonically equivalent (see Subsection 1.7).

The marked vector $\Xi$ is $K$-invariant,
therefore (see Subsection 1.6) we can consider
the matrix element
$$\b_\alpha(g):=<T_\alpha(g)\Xi,\Xi>$$
as a function on $G/K$ or a function on $K\setminus G/K$.
Vectors $\Xi$ and $T_\alpha(g)\Xi$ are elements
of the supercomplete basis in $H[M_\alpha]$,
 therefore the function
$\b_\alpha$ can be easily evaluated.

In the matrix ball model of $G/K$ the function
$\b_\alpha$ is given by the formula
\begin{equation}
\b_\alpha(z)=\det(1-zz^*)^{\alpha/2}; \qquad z\in\B_{p,q}
\end{equation}

Let us obtain the formula for $\b_\alpha$ as function on $K\setminus G/K$.
Denote by $a_t$ the element of $\O(p,q)$ given by the matrix
\begin{equation}
a_t=
\begin{pmatrix}
\begin{array}{ccc}\ch t_1&  & \\ &\ddots&\\ && \ch t_p \end{array}&
\begin{array}{ccc}\sh t_1&  & \\ &\ddots \\ && \sh t_p\end{array}&
\begin{array}{ccc}0&\cdots&0\\ & \ddots & \\0& \cdots& 0 \end{array}\\
\begin{array}{ccc}\sh t_1&  & \\ &\ddots \\ && \sh t_p\end{array}&
\begin{array}{ccc}\ch t_1&  & \\ &\ddots&\\ && \ch t_p \end{array}&
\begin{array}{ccc}0&\cdots&0\\ & \ddots & \\0& \cdots& 0 \end{array}\\
\begin{array}{ccc}0\quad&\cdots&0\\ & \ddots & \\0\quad& \cdots& 0 \end{array}&
\begin{array}{ccc}0\quad&\cdots&0\\ & \ddots & \\0\quad& \cdots& 0 \end{array}&
\begin{array}{ccc}1&\cdots&\\ & \ddots & \\ & \cdots& 1 \end{array}
\end{pmatrix}
\end{equation}
It is easy to show that arbitrary element $g$ of $G=\O(p,q)$
can be represented in the form
$$g=k_1a_t k_2;\qquad \mbox{where}\qquad k_1,k_2\in K$$
The collection of parameters $t=(t_1,\dots,t_p)$
is uniquely defined up to permutations of $t_j$ and
reflections
\begin{equation}
(t_1,\dots,t_p)\mapsto (\sigma_1 t_1,\dots,\sigma_pt_p)
\end{equation}
where $\sigma_j=\pm1$.

We denote by $\cal A$ the subgroup in $\O(p,q)$ consisting of all elements
$a_t$. We denote by $D_p$ the group of  transformations of $\R^p$
generated by permutations of coordinates and reflections (1.21).

We identify the set $K\setminus G/K$
with the set of $D_p$-orbits on $\cal A$.

 In coordinates $(t_1,\dots,t_p)$ the matrix element $\b_\alpha$
is given by the formula
\begin{equation}
\b_\alpha(t_1,\dots,t_p)=\prod_{k=1}^p\ch^{-\alpha}t_k
\end{equation}
Hence, we must obtain the expansion (1.18) of the function $\b_\alpha$
given by  formula (1.19) or (1.22)
 in positive definite spherical
functions.

 Our  purpose in Section C is to give an expression
for spherical functions.

\bol

{\bfff C. Spherical representations and spherical transform}

\nopagebreak

\bol

{\bff 1.\punkt. Parabolic subgroup.} Consider the
space $\R^p\oplus \R^q$ equipped with the indefinite symmetric
form $J$ defined by  formula (1.5).
 A subspace $V\subset \R^p\oplus \R^q$ is called
{\it isotropic} if the form $J$ is zero on $V$.

An {\it isotropic flag} $\cal V$ in $\R^p\oplus \R^q$
is a family of isotropic subspaces
$${\cal V}:\quad V_1\subset V_2\subset\dots \subset V_p;
\qquad \mbox{where}\quad \dim V_j=j$$
The {\it flag manifold}  $\cal F$ is the space of all isotropic flags
 in $\R^p\oplus \R^q$.

The space $\cal F$ is an $\O(p,q)$-homogeneous space.
A {\it minimal parabolic subgroup}
is the stabilizer of a point in $\cal F$.
Let us give more explicit description of
the minimal parabolic subgroup.

For this let us consider the basis
$v_1,\dots, v_p,w_1,\dots,w_{q-p},
v_1',\dots,v_p'$
 in $\R^p\oplus\R^q$ defined by
\begin{equation}
v_j=\frac 1{\sqrt{2}}(e_j+e_{q+j});\qquad
v_j'=\frac 1{\sqrt{2}}(e_j-e_{q+j});\quad w_k=e_{p+k}
\end{equation}
Then
$$
J(e_k,e_k')=1; \qquad J(f_k,f_k)=1
$$
and the scalar products of all other pairs of
basic vectors  are zero.

Denote by ${ L}_k$ the subspace in $\R^p\oplus \R^q$
generated by the basic vectors $v_1,\dots,v_k$.
We denote  by $P\subset \O(p,q)$
the stabilizer in $\O(p,q)$ of the isotropic flag
\begin{equation}{\cal L}:\,\,\,{L}_1\subset\dots\subset{ L}_p
\end{equation}
The subgroup $P$ is a minimal parabolic   subgroup in $\O(p,q)$
and
$${\cal F}\simeq \O(p,q)/P$$

 Elements of      the
parabolic subgroup $P$ in the basis (1.23) have the form
\begin{equation}
\begin{pmatrix}A&*&*\\0&C&*\\0&0&A^{t-1}\end{pmatrix}
\end{equation}
where $A$ is an upper triangular matrix and $C\in\O(q-p)$.

Elements of the subgroup $\cal A$
(see Subsection 1.15) in new basis are diagonal
matrices with eigenvalues
$$
e^{t_1},\dots, e^{t_p},1,\dots, 1,e^{-t_1},\dots, e^{-t_p}
$$

\mal

{\scc Remark.} Let us change the order of basic elements
(1.23) to
$v_1,\dots, v_p$ ,$w_1,\dots,w_{q-p}$,
$v_p',\dots,v_1'$.
Then elements of the parabolic subgroup $P$ will be upper triangular matrices.

    \mal

{\bff 1.\punkt. Spherical representations.}
Denote by $\Gr_k$ the space of all $k$-dimensional
isotropic subspaces in $\R^p\oplus\R^q$.
Consider the tautological embedding
of the flag space $\cal F$ to the product of the Grassmannians
$\times_{k=1}^p\Gr_k$ ( to each point
${\cal V}:V_1\subset\dots\subset V_p$ we assign the point
$(V_1,\dots,V_p)\in\times_{k=1}^p\Gr_k$).

Consider the natural action
of $\O(p,q)$ on $\Gr_k$. For $g\in G$  we denote by $j_k(g,V)$
the Jacobian of the transformation $g$ at the point $V\in \Gr_k$.
By $J(g,{\cal V})$ we denote the Jacobian of the transformation
$g$ on the flag space $\cal F$.

Fix $s_1,\dots,s_p\in \C$. Assume $s_0=s_{p+1}=0$.
 We define the representation
$\widetilde\pi_s$ of the group $\O(p,q)$
 in the space of functions  on $\cal F$
by the formula
$$\widetilde\pi_s(g)f(V_1,\dots,V_p)=f(gV_1,\dots,gV_p) J(g,{\cal V})^{1/2}
  \prod_{k=1}^p j_k(g,V_k)^{(s_{j-1}-2s_j+s_{j+1})/2}
$$

 {\scc Remark.} The representation
$\widetilde\pi_s$ is a Harish-Chandra
module%
\footnote{This means
(for instance, see \cite{Kna}) that the spectrum of the maximal
compact subgroup $K$ in the space of functions
on $\cal F$ has finite multiplicities.
 } and hence a topology
in the space of functions on $\cal F$ is not essential.
For instance,  we can consider the space $L^2({\cal F})$,
the space of smooth
functions $C^\infty({\cal F})$,
 the space of distributions ${\cal D}({\cal F})$,
 the space of hyperfunctions etc.

\mal

  {\scc Remark.} Consider the $\delta$-function $\delta_{\cal L}$
supported at the point ${\cal L}\in\cal F$ (see (1.24)).
 It is easy to observe that
the function  $\delta_{\cal L}$ is an eigenfunction of $P$ and for
a matrix $g\in P$ given by (1.25) we have
\begin{equation}
\widetilde\pi_s(g)\delta_{\cal L}({\cal V})=
 \exp\bigl\{\sum_{j=1}^p t_j (s_j-(q+p)/2+j)\bigr\}
           \delta_{\cal L}({\cal V})
\end{equation}
where $e^{t_j}$ are the absolute values of
the eigenvalues of the matrix $A$ (see (1.25)).

\mal

We want to define a canonical irreducible subquotient $\pi_s$
in  $\widetilde\pi_s$.

\mal

{\scc Remark.} For generic $s\in\C^n$ the representation
 $\widetilde\pi_s$ is irreducible
and hence $\pi_s\simeq\widetilde\pi_s$.

\mal

Consider the function $f_0({\cal V})=1$ on the space $\cal F$,
it is the unique $K$-invariant function on $\cal F$
(since $\cal F$ is $K$-homogeneous). Denote by
$S$ the cyclic span of $f_0$. Denote by $R$ the sum of all proper
$\O(p,q)$-submodules in $S$.

\mal

{\scc Lemma 1.\fact.} $R\ne S$.

{\scc Proof.}
Indeed there is the unique $K$-fixed vector in $S$ and this vector
is cyclic. Hence, it can't be element of a proper submodule.
Hence, a proper submodule in $S$ hasn't $K$-fixed vector.
Hence, $R$ also hasn't $K$-fixed vectors and hence
$f_0\notin R$.\hfill$\boxtimes$

  We define  the $\O(p,q)$-module $\pi_s$ by
$$\pi_s=S/R$$

{\scc Theorem 1.\fact.}%
\footnote{for instance, see \cite{Hel}, Theorem 4.4.3}
{\it The representations $\pi_s$ are precisely all spherical representations
of $\O(p,q)$. Moreover}
$$\pi_s\simeq\pi_{s'} \quad \mbox{iff there exists}\quad  \gamma\in D_p
\quad\mbox{  such that}\quad \gamma s=s'$$

Hence, we can consider our Plancherel measure $\nu_\alpha$ as
a measure on $\C^p/D_p$. {\it It will be more convenient for us
to consider the Plancherel measure as a $D_p$-invariant
measure on $\C^p$ or any measure on on $\C^p$ whose $D_p$-average
is $\nu_\alpha$.}

\mal

{\bff 1.\punkt. Unitary spherical representations.}
{\scc Lemma 1.\fact.} {\it Assume the representation $\pi_s$ be  unitary.
Then for any $j$}
\begin{equation}
\Re s_j=0 \qquad \mbox{or} \qquad \Im s_j=0
\end{equation}

{\scc Proof.} The representation dual to $\pi_s$ is $\pi_{-s}$.
The complex conjugate representation to $\pi_s$ is $\pi_{\ov s}$.
If $\pi_s$ is unitary, then the dual representation is equivalent to the
complex conjugate representation. Hence,
$-s=\gamma \ov s$ for some $\gamma\in D_p$.\hfill$\boxtimes$

\mal

 If $s_1,\dots,s_p$ are pure imaginary,
then the representation $\pi_s$ is  unitary
in $L^2({\cal F})$.  These representations are called
{\it representations of the principal nondegenerate series}.

For some other values of $s$ representations $\pi_s$
also are unitary, but scalar product in these cases is more complicated.

{\scc Theorem 1.\fact.} ({\it see} \cite{Hel}, 4.8.1)
{\it Denote by $\rho$ the vector
$$\bigl((q+p)/2-1,(q+p)/2-2,\dots, (q-p)/2)\bigr)\in\R^p$$
Denote by $Q$ the convex polyhedron in $\R^p$ with vertices
$\gamma\rho$ where $\gamma \in D_p$. Then for each unitary
representation $\pi_s$                         }
\begin{equation}
(\Re s_1,\dots,\Re s_p)\in Q
\end{equation}
{\it Moreover, spherical function of a spherical representation
$\pi_s$ is bounded if and only if
condition (1.28) holds.}

Our next purpose is to give the integral formula
for the spherical functions in an explicit form.
For this we must give another realization
of $G/K$.

{\bff 1.\punkt. Matrix wedges.} First, consider the case $p=q$.
Consider the matrix ball $\B_{q,q}(\R)$. Consider the {\it Cayley transform}
\begin{equation}
\Cay:z\mapsto\frac{1-z}{1+z}
\end{equation}
Then the map $\Cay$ transfers the matrix ball $\B_{q,q}(\R)$ to the wedge
$\W_q$ consisting of matrices $R$ satisfying the condition%
\footnote{A matrix satisfying the condition $R+R^t>0$ is called
{\it dissipative.}}
$$ \qquad R+R^t>0$$
(where the notation $Q>0$ means that a matrix $Q$ is
positive definite).
It is convenient to wright $R$ in the form
$$R=T+S\qquad \mbox{where}\qquad T=T^t>0;\quad S=-S^t$$

The group $\O(q,q)$ acts on $\B_{q,q}(R)$ and hence it acts
on $\W_{q,q}$. For description of the last action
we
consider the basis (1.23) in $\R^q\oplus\R^q$. In our case $p=q$,
and hence the basic elements $w_j$ are lacking. Hence,
  $\O(q,q)$ becomes  the group of
real $(q+q)\times (q+q)$-matrices having the form
$g=\begin{pmatrix}a&b\\c&d\end{pmatrix}$ and  satisfying the condition
$$\begin{pmatrix}a&b\\c&d\end{pmatrix}
\begin{pmatrix}0&1\\1&0\end{pmatrix}
\begin{pmatrix}a&b\\c&d\end{pmatrix}^t=
\begin{pmatrix}0&1\\1&0\end{pmatrix}
$$
The group $\O(q,q)$ acts on $\W_q$ by fractional linear
transformations
$$R\mapsto R^{[g]}:= (a+Rc)^{-1}(b+Rd)$$
In this model, the parabolic subgroup $P\subset \O(q,q)$
becomes the group of real matrices having the form
\begin{equation}
\begin{pmatrix}a&b\\0&a^{t-1}\end{pmatrix};\qquad\mbox{where}\qquad d
\,\,\, \mbox{is upper triangular}
\end{equation}
Hence, the parabolic subgroup acts on $W_q$ by affine transformations
$$R\mapsto a^{-1}Ra^{t}+a^{-1}b$$
We emphasis that $a^{-1}b$ is a skew-symmetric matrix.

We can easily wright the eigenfunctions of the group $P$ on the wedge $\W_q$
\begin{equation}
\Psi_{s_1,\dots,s_q}(R)=\prod_{j=1}^q \det[T]_j^{(-\theta_j+s_j-s_{j+1})/2}
\quad\mbox{where}\quad \theta_1=\dots=\theta_{p-1}=1,\,\, \theta_p=0
\end{equation}
Here the symbol $[T]_j$ denotes the left upper $j\times j$ block of the matrix
$T$.

Consider  $g\in P$  given by  formula (1.30). Let $e^{t_1},\dots,
e^{t_p}$ be the absolute values of the diagonal elements of the block $a$. Then
\begin{equation}
\Psi_{s_1,\dots,s_q}(R^{[g]})
=\exp\bigl\{\sum_{j=1}^q(s_j-q+j)t_j\bigr\} \Psi_{s_1,\dots,s_q}(R)
\end{equation}

{\scc Remark.} Compare (1.32) and (1.26).

\mal

{\bff 1.\punkt. Sections of wedges.}
 Consider arbitrary group $\O(p,q)$. Let us represent
a point $z\in\B_{p,q}(\R)$ as  block  $p\times(p+(q-p))$-matrix
$z=(z_1 \,\, z_2)$. Consider  block $(p+(q-p))\times(p+(q-p))$ matrix
$$\widetilde z =\begin{pmatrix}0&0\\z_1 &z_2\end{pmatrix}\in \B_{q,q}(\R)$$
Thus, we realized the matrix ball $\B_{p,q}(\R)$ as a submanifold of
$\B_{q,q}(\R)$. The image $\SW_{p,q}$ of $\B_{p,q}(\R)$
 under the Cayley transform (1.29)
is the set of $(p+(q-p))\times(p+(q-p))$-matrices
 $R\in \W_{q}$ having block structure
\begin{equation}
R=\begin{pmatrix}1&0\\Q&H\end{pmatrix}
\end{equation}
The condition $R+R^t>0$ for matrix (1.33) is equivalent to the condition
\begin{equation}
\frac12(H+H^t)-LL^t>0
\end{equation}
(spaces $\SW_{p,q}$ are real sections of so-called
 Siegel domains of the second type, see \cite{Pya})

We will wright matrices $R\in \SW_{p,q}$ in the form
$$R=\begin{pmatrix}1&0\\2L& M+N\end{pmatrix};\qquad  M=M^t, N=-N^t$$

The condition (1.34) can be represented in the form
\begin{equation}
M-LL^t>0 \qquad \mbox{or} \qquad
\begin{pmatrix}1&L^t\\L&M\end{pmatrix}>0
\end{equation}

The eigenfunctions of the
 parabolic subgroup $P\subset\O(p,q)$ in this model
are given by the formula
\begin{equation}
\Psi_{s_1,\dots,s_p}(R)=\prod_{j=1}^p
\det\Bigl[\begin{pmatrix}1&L^t\\L&M\end{pmatrix}
     \Bigr]_{q-p+j}^{(-\theta_j+s_j-s_{j+1})/2}
=\prod_{j=1}^p \det[M-LL^t]_j^{(-\theta_j+s_j-s_{j+1})/2}
\end{equation}
where $\theta_1=\dots=\theta_{p-1}=1,\,\, \theta_p=\half(q-p)$
(see simple calculations in \cite{Ner4})

The Berezin kernel $L_\alpha$ (see Subsection 1.10)
 in the models $\W_q$, $\SW_{p,q}$
is given by the formula
$$L_\alpha(R_1,R_2)=
\Bigl[\frac{\det(R_1+R_2^t)}
    {\det(1+R_1)\det(1+R_2)}\Bigr]^{-\alpha}
$$
This gives the following expression for the function $\b_\alpha(R)$
in our coordinates
\begin{equation}
\b_\alpha(R)=\frac{\det
2^{\alpha}\begin{pmatrix}1&L^t\\L&M\end{pmatrix}^{\alpha/2}}
{\det(1+M+N)^{\alpha}}
\end{equation}

The $\O(p,q)$-invariant measure on $\SW_{p,q}$ is
$$\det\begin{pmatrix}1&L^t\\L&M\end{pmatrix}^{-(p+q)/2}dL\,dM\,dN$$
where
$dL$, $dM$, $dN$ are Lebesgue measures on the spaces of matrices.

\mal

{\bff 1.\punkt. Canonical embedding of the spherical
 $G$-module $\pi_s$ to the space
$C^\infty(G/K)$.}
First, we define the canonical intertwining
operator $J_s$ from $\widetilde\pi_s$ to $C^\infty(G/K)$.
This operator is uniquely defined by the following property
$$J_s:\quad \delta_{\cal L}\mapsto \Psi_s$$
where $P$-eigenfunctions  $\delta_{\cal L}$, $\Psi_s$
were defined in Subsections 1.17, 1.19--1.20.
 By the intertwining property  we obtain
$$J_s\delta_{g{\cal L}}(R)=\Psi_s(R^{[g]})$$
and this defines the operator $J_s$ on all
$\delta$-functions. Then we extend $J_s$ by linearity and continuity
to the whole space of distributions on $\cal F$.

\mal

{\scc Lemma 1.\fact.} {\it The operator $J_s$ induces an embedding
of the subquotient $\pi_s$ to $C^{\infty}(G/K)$.}

\mal

Le us denote by
$d{\bf k}$ the Haar measure on $K=\O(p)\times\O(q)$.
  We assume that the measure of the whole group is 1.

\mal

{\scc Proof.} Let $R$, $S$ be the same as in Subsection 1.17.
 Let $Q\subset C^{\infty}(G/K)$ be a $G$-invariant
closed subspace. Then for any function $f\in R$, its {\it average}
$$f^K(R)=\int_{{\bf k}\in K} f([R]^{\bf k})d{\bf k}$$
  is contained in $Q$. Hence, $Q$ contains
a $K$-invariant function.

  By this reason, $J_s$ maps the submodule $R$ to 0
(since $R$ hasn't $K$-invariants). Assume that $J_s$ is
zero on $S$. Then $J_s$ is an operator from $\widetilde\pi_s/S$
to $C^{\infty}(G/K)$. But  the module $\widetilde\pi_s/S$
hasn't $K$-invariant vectors. Hence, $J_s$ is  identical zero
and this contradicts to its definition. \hfill$\boxtimes$

\mal

Obviously the $K$-fixed function $f_0=1$ on $\cal F$
 can be represented in the form
$$f_0({\cal V})=\int_{{\bf k}\in K}
  \delta_{\cal L}({\bf k}{\cal V})d{\bf k}$$
 Hence, its image under $J_s$ is the $K$-average of
$\Psi_s$. This gives the integral formula
for spherical function given in the next Subsection.

{\bff 1.\punkt. Integral formula for spherical functions.}
Spherical functions are $K$-averages of  $P$-eigenfunctions
on $G/K$
\begin{equation}
\Phi_{s_1,\dots,s_p}(R)=\int_{{\bf k}\in\O(p)\times\O(q)}
  \Psi_{s_1,\dots,s_p}(R^{[{\bf k}]})d\pi({\bf k})
\end{equation}

{\scc Lemma 1.\fact.}
\begin{equation}
|\Phi_{s_1,\dots,s_p}(t)|\le
\Phi_{\Re s_1,\dots,\Re s_p}(t)
\end{equation}

{\scc Proof} is obvious. \hfill $\boxtimes$

{\bff 1.\punkt. Spherical transform.}
Let $f(z)$ be a $K$-invariant function on $G/K$. Then the {\it spherical
transform} of $f$ is defined by the formula
\begin{equation}
\widehat f(s)=\int_{G/K}\Phi_s(z)f(z)d\lambda(z)
\end{equation}
where $\lambda$ is the $G$-invariant measure   on $G/K$.

If $f\in L^2\cap L^1(G/K)$, then the
 {\it Gindikin--Karpelevich inversion formula}
(see \cite{GK1}, \cite{GK2}, \cite{Hel}, \cite{FK})
is valid
\begin{equation}
f(z)= C\cdot\int_{i\R^p}\widehat f(s)\Phi_s(z)\r(s)ds
\end{equation}
where $C$ is a known constant (see \cite{Hel}, formula (4.6.40))
and $\r(s)$ is the {Gindikin--Karpelevich density}. For $G=\O(p,q)$
it is given by the formula
\begin{align}
&\cdot\r(s)=
 \prod\limits_{k=1}^p\frac{\Gamma((q-p)/2+s_k) \Gamma((q-p)/2-s_k)}
                       {\Gamma(s_k)\Gamma(-s_k)} \times\\ \times
&\prod\limits_{1\le k<l\le p}\frac
        {\Gamma(\frac12(1+s_l+s_k))\Gamma(\frac12(1+s_l-s_k))
        \Gamma(\frac12(1-s_l+s_k))\Gamma(\frac12(1-s_l-s_k))}
 {\Gamma(\frac12(s_l+s_k))\Gamma(\frac12(s_l-s_k))
        \Gamma(\frac12(-s_l+s_k))\Gamma(\frac12(-s_l-s_k))}
\end{align}

{\scc Remark.} This expression is an elementary function.
For instance, using the complement formula for $\Gamma$, we
reduce factor (1.43)  to the form
\begin{equation}\prod_{1\le k<l\le p}
  (s_k^2-s_l^2)\tg \pi(s_k+s_l)\tg \pi(s_k-s_l)
\end{equation}
If $q-p$ is even, then (1.42) equals to
\begin{equation}
\prod_{k=1}^p \Bigl\{\prod_{\tau=0}^{(q-p)/2-1}(\tau^2-s_k^2)\Bigr\}
\end{equation}
If $(q-p)$ is odd, then (1.42) equals to
\begin{equation}
\prod_{k=1}^p \Bigl\{
s_k\tg \pi s_k
\prod_{\tau=0}^{(q-p-3)/2}
  ((\tau+1/2)^2-s_j^2)\Bigr\}
\end{equation}

For pure imaginary $s$ we can replace (1.42)-(1.43) by
$$\prod_{k=1}^p
\Bigl|\frac{\Gamma((q-p)/2+s_k)}{\Gamma(is_k)}\Bigr|^2
\prod\limits_{1\le k<l\le p} \Bigl| \frac
        {\Gamma(\frac12(1+s_l+s_k))\Gamma(\frac12(1+s_l-s_k))}
 {\Gamma(\frac12(s_l+s_k))\Gamma(\frac12(s_l-s_k))}\Bigr|^2
$$
since $s_k$ is imaginary. Nevertheless long expression
(1.42)--(1.43) is more convenient for our calculations.

               \mal

{\bff 1.\punkt. Another formula for spherical transform.}
By   integral formula for spherical functions
(1.38) we can wright  spherical transform (1.40) in the
following  form
\begin{equation}
\widehat f(s)=\int_{G/K}f(z)\Psi_s(z)d\lambda(z)
\end{equation}

\mal

{\bff 1.\punkt. Further structure of the paper.}
We want to obtain an expansion of the function $\b_\alpha(z)$
in spherical functions. If $\b_\alpha\in L^1\cap L^2(\B_{p,q}(\R))$
(or $\alpha>p+q-1$),
  then it is sufficient to evaluate
the spherical transform of the function
$\b_\alpha(z)$, and the Gindikin-Karpelevich inversion
formula gives required expansion.

 In  Section 3 we  evaluate the spherical transform of     $\b_\alpha(z)$
using  formula (1.47).
The final result is given in Theorem 2.1.
Then in \S4 we  construct the analytic continuation of
 our formula to arbitrary $\alpha$.  As result, we obtain
an expansion of $\b_\alpha$ in spherical functions.
 In \S 5 we prove positive
definiteness of these spherical functions.

\bol

{\bfff D. Deformation of $L^2$ on Riemannian compact symmetric space
and kernel representations of $\O(p+q)$.}

\nopagebreak

\bol

This subject is a supplement to the main topic of the paper.

\mal

{\bff 1.\punkt. The symmetric spaces $\U(p+q)/\U(p)\times\U(q)$
and $\O(p+q)/\O(p)\times\O(q)$.}
Consider the group $\U(p+q)$ consisting of all complex block
$(p+q)\times (p+q)$-matrices
$\begin{pmatrix}a&b\\c&d\end{pmatrix}$
satisfying the condition
$$\begin{pmatrix}a&b\\c&d\end{pmatrix}
\begin{pmatrix}1&0\\0&1\end{pmatrix}
\begin{pmatrix}a&b\\c&d\end{pmatrix}^*=
\begin{pmatrix}1&0\\0&1\end{pmatrix}$$
 Consider the subgroup $\O(p+q)\subset\U(p+q)$
consisting of real matrices.

Consider the Grassmannians $\Gr_{p,q}(\C)$ and
$\Gr_{p,q}(\R)$ consisting of
$p$-dimensional subspaces in $\C^{p+q}$ and
$\R^{p+q}$ respectively. Obviously, these
 Grassmannians are symmetric spaces
\begin{gather*}
\Gr_{p,q}(\C)=\U(p+q)/\U(p)\times\U(q)\\
\Gr_{p,q}(\R)=\O(p+q)/\O(p)\times\O(q)
\end{gather*}

Denote by $\M_{p,q}(\C)$ (resp. $\M_{p,q}(\R)$)
 the space of all $p\times q$-matrices over $\C$ (resp. over $\R$).
For any $z\in \M_{p,q}$ we define its graph
 $\graph_z\subset \Gr_{p,q}$.
Obviously the map $z\mapsto\graph_z$ is
an embedding of $\M_{p,q}\to\Gr_{p,q}$ and the image
of the embedding is dense in $\Gr_{p,q}$.

In the coordinate $z\in\M_{p,q}$,
 the action of the group $\U(p+q)$ on Grassmannian
  is given by
the formula
\begin{equation}
z\mapsto z^{[g]}=(a+zc)^{-1}(b+zd)
\end{equation}
 coinciding with  formula (1.7).

    \mal

{\bff 1.\punkt. Representations $\widetilde T_{-n}$.}
Fix $n=0,1,2,\dots$.
Denote by $\phi_a(z)$ the polynomial on $\M_{p,q}$
given by the formula
$$\phi_a(z)=\det(1+a^*z)^n; \qquad\mbox{where}\quad a\in \M_{p,q}$$
Denote by $H_{-n}$ the linear span of
 all polynomials $\phi_a(z)$.
Obviously the space $H_{-n}$ is finite dimensional
(since degree of the polynomial $\phi_a(z)$ is $pn$).

Consider the action of the group $\U(p,q)$
in the space $H_{-n}$
given by the formula
\begin{equation}
\widetilde T_{-n}\begin{pmatrix}a&b\\c&d\end{pmatrix}f(z)=
f\bigl((a+zc)^{-1}(b+zd)\bigr)\det(a+zc)^n
\end{equation}
coinciding with  formula (1.10).
It is easy to check, that the transformations $\widetilde T_{-n}(g)$
preserve the space $H_{-n}$.

Consider the scalar product in $H_{-n}$
given by

$$<f_1(z),f_2(z)>_{-n}=
C_n\int_{\M_{p,q}}f_1(z)\ov {f_2(z)}\det(1+z^*z)^{-n-p-q}dz$$
where the normalization constant $C_n$ is defined by the condition
$<1,1>_{-n}=1$.

It is easy to check that the operators $\widetilde T_{-n}(g)$
are unitary with respect to this scalar product and
\begin{equation}
<\phi_a,\phi_b>_{-n}=\det(1+a^*b)^n
\end{equation}
Hence, the (finite dimensional) hilbert space $H_{-n}$ is
the hilbert space $H^\circ$ defined by
the positive definite kernel
\begin{equation}
L_{-n}(a,b)=\det(1+a^*b)^n  ;\qquad a,b\in\M_{p,q}
\end{equation}

{\bff 1.\punkt. Kernel representations of $\O(p+q)$.}
The {\it kernel representation} $T_{-n}$
of the group $\O(p+q)$ is the restriction of the
representation $\widetilde T_{-n}$ to the subgroup
$\O(p+q)$.

The also define the {\it marked vector} $\Xi$
$$\Xi:\qquad f(z)=1$$

\mal

{\bff 1.\punkt. Limit as $n\to \infty$.} Let us consider
the kernel
$$M_{-n}(z,u)=\frac{\det(1+z^tu)^n}{\det(1+z^tz)^{n/2}\det(1+u^tu)^{n/2}}$$
on $\M_{p,q}(\R)$. By Lemma 1.5.e) the kernel $M_{-n}$ is positive definite.
A simple calculation shows,
 that  the kernel is $\O(p+q)$-invariant.
Consider the hilbert space $H^\circ[M_{-n}]$. The operator
$$Af(z)=\det(1+z^*z)^{n/2}f(z)$$
defines the canonical unitary $\O(p+q)$-intertwining operator
$H^\circ[L_\alpha]\to H^\circ[M_\alpha]$.

 The
  arguments  given in Subsection 1.13 show
that a natural limit of the spaces $H_{-n}$ as $n\to\infty$
is
$$L^2\bigl(\O(p+q)/\O(p)\times\O(q)\bigr)$$

\mal

{\bff 1.\punkt. Preliminary remarks on the Plancherel formula.}
By Subsection 1.7 the matrix element
$$\b_{-n}(g)=<T_{-n}(g) \Xi,\Xi>_{H_{-n}}$$
is a function on $\O(p+q)/\O(p)\times\O(q)\simeq\Gr_{p,q}(\R)$.
In the coordinate $z\in \M_{p,q}(\R)$ it is given by
$$\b_{-n}(z)=\det(1+zz^t)^{-n/2}$$

Denote by $\widehat{\O(p+q)}_{sph}$ the set of all irreducible
representations of $\O(p+q)$ having
an $\O(p)\times\O(q)$-invariant vector ({\it spherical vector});
description of this (countable) set is given by Helgason theorem,
\cite{Hel}, Theorem 5.4.1. By $H_\rho$ we denote
the space of a spherical representation
$\rho\in \widehat{\O(p+q)}_{sph}$. Denote by
$\xi_\rho$ the spherical vector in $H_\rho$ having
unit length.

Arguments given in Subsections 1.13-1.14 show, that the decomposition
of $T_{-n}$ in irreducible representations has the form
$$T_{-n}(g)=\bigoplus_{\rho\in\Delta_n} \rho(g)$$
where $\Delta_n$ is a finite subset in $\widehat{\O(p+q)}_{sph}$.

The scalar product in $\bigoplus_{\rho\in\Delta_n} H_\rho $
has the  form
\begin{equation}<\bigoplus_{\rho\in\Delta_n} v_\rho,
\bigoplus_{\rho\in\Delta_n} w_\rho>=\sum_{\rho\in\Delta_n}
\nu^n_\rho <v_\rho,w_\rho>_{H_\rho}
\end{equation}
where $v_\rho,w_\rho\in H_\rho$ and $\nu^n_\rho$ are positive constants.
The formula (1.52) is called the {\it Plancherel formula}.

We normalize the constants $\nu^n_\rho$ by the assumption
$$ \mbox{the image of} \quad \Xi \quad \mbox{in}\quad
 \bigoplus_{\rho\in\Delta_n} H_\rho  \quad \mbox{is}\quad
 \bigoplus_{\rho\in\Delta_n} \xi_\rho$$
The constants $\nu^n_\rho$ are  evaluated in Section 2 as a corollary
of the Plancherel formula for kernel-representations
of $\O(p,q)$.

\bol

{\bfff E.
 An interpolation between
$$\bf L^2(\O(p,q))/\O(p)\times \O(q)) \,\,\, \mbox{\large\bf and} \,\,\,
L^2(\O(p+q))/\O(p)\times \O(q))\,\,\,\mbox{\large\bf?}$$}

\nopagebreak

\mal

The purpose of the Section is a formulation of a strange problem.

\mal

{\bff 1.\punkt. General representations $\widetilde T_\alpha$.}
Denote by $Hol(\B_{p,q})$ the space of holomorphic functions
in $\B_{p,q}(\C)$ equipped with topology of the uniform convergence on
compacts.

Consider  arbitrary $\alpha\in\C$ and consider  the
action of $\U(p,q)$ in $Hol(\B_{p,q})$  given by the formula
\begin{equation}
\widehat T_{\alpha}\begin{pmatrix}a&b\\c&d\end{pmatrix}f(z)=
f\bigl((a+zc)^{-1}(b+zd)\bigr)\det(a+zc)^{-\alpha}
\end{equation}
Denote by ${\cal H}^\circ_\alpha$ the cyclic span of the function
$f(z)=1$. Denote by $\widetilde T_{\alpha}$ the restriction of
$\widehat T_{\alpha}$ to  ${\cal H}^\circ_\alpha$.
It is easy to observe  that  $\widetilde T_{\alpha}$
is an irreducible Harish-Chandra module.

Let us denote by $H^{fin}_\alpha$ the space of polynomials contained
in ${\cal H}^\circ_\alpha$. Consider the action of the Lie algebra
$u(p,q)$ in  $H^{fin}_\alpha$. The space $H^{fin}_\alpha$
is an irreducible $u(p,q)$-module with a highest weight.
For $\alpha\in\R$ there exists the unique $u(p,q)$-invariant
hermitian form in $H^{fin}_\alpha$ (it is called  {\it Shapovalov form}%
\footnote{Its definition for highest weight modules
of various group is uniform,
see for instance \cite{Ner2}.}).
In general this form is indefinite.

If $\alpha\in\R$ satisfies  Berezin conditions (1.9), then the
Shapovalov form is positive definite. It
 coincides with the Berezin scalar product,
and representations $\widetilde T_\alpha$ coincides with
 representations   $\widetilde T_\alpha$ constructed in Subsection 1.10.
 If $\alpha$
is a negative integer, then $\widetilde T_\alpha$
is finite dimensional. By the unitary Weyl
trick, there is no difference
between finite dimensional representations of $\U(p,q)$,
holomorphic finite dimensional representations of $\GL(p+q,\C)$
and finite dimensional representations of $\U(p+q)$.
The representations $\widetilde T_\alpha$  for negative
integer $\alpha$ differs from the representations $\widetilde T_{-n}$
from Subsection 1.27  by a nonessential change of notations.

\mal

{\bff 1.\punkt. Nonunitary kernel representations of $\O(p,q)$?}
Consider the  restriction  $T_\alpha$ of  $\widetilde T_\alpha$
to the subgroup $\O(p,q)$. It is a well-defined representation
of the group $\O(p,q)$ in the space ${\cal H}^\circ_\alpha$.

We  have seen that
\begin{align*}
\lim_{\alpha\to+\infty} T_\alpha&\simeq
           L^2\bigl(\O(p,q))/\O(p)\times \O(q)\bigl)\\
\lim_{n\to-\infty} T_n&\simeq
L^2\bigl(\O(p+q))/\O(p)\times \O(q)\bigl)
\end{align*}

It seems that the Plancherel formula (2.5)-(2.15)
gives the decomposition of the kernel representation $T_\alpha$
for any complex $\alpha$. Unfortunately
it is a result of 'mathematical physics level' . This is the solution of
a problem which hasn't a satisfactory formulation (since
the definition of the abstract Plancherel formula
doesn't exist for nonunitary representations).

\mal

{\scc Remark.} For the case $p=1$ the space ${\cal H}^\circ_\alpha$
equipped with the Shapovalov form is a Pontryagin space%
\footnote{This means that negative inertia index of the Shapovalov form
is finite.}
and in this case our Plancherel formula is really the Plancherel
formula for arbitrary real $\alpha$.

                            \mal

The questions of this type are discussed for a long
time and they  arise
to deep Molchanov work \cite{Mol}(1980) containing
the Plancherel decomposition of tensor products of
unitary representations of ${\rm SL}(2,\R)$ (see also \cite{MD1}
and references in this paper).
 It was clear, that Molchanov formula  give
formal interpolation between tensor products of unitary
representations of ${\rm SL}(2,\R)$ and tensor products of
finite dimensional
representations.
Unfortunately before our time  a quite satisfactory
group-theoretical interpretation of this interpolation
 doesn't exist.

\bol

{\bfff  2. Formulation of results.}

\nopagebreak

\bol

\addtocounter{sec}{1}
\setcounter{equation}{0}
\setcounter{fact}{0}
\setcounter{punkt}{0}

{\bff  2.\punkt. Large $\alpha$.}
{\sc Theorem 2.\fact.}{\it Let $\alpha>(p+q)/2-1$. Then
the spectrum of the kernel-representation
 $T_\alpha$ is supported by nondegenerate principal
series and the Plancherel decomposition is given by the formula
\begin{align}
&\prod_{k=1}^p\ch^{-\alpha} t_j=
C\cdot{2^\alpha}\frac1{\prod_{j=1}^p{\Gamma(\alpha-j+1)}}\times\\&
\int_{i\R^p}\prod_{k=1}^p \left\{\Gamma(\frac12(\alpha-(p+q)/2+1+s_k))
              \Gamma(\frac12(\alpha-(p+q)/2+1-s_k))\right\}\times\\&
\times     \prod_{k=1}^p\frac{\Gamma((q-p)/2+s_k) \Gamma((q-p)/2-s_k)}
                       {\Gamma(s_k)\Gamma(-s_k)}  \times \\&
\prod_{1\le k<l\le p}\frac
        {\Gamma(\frac12(1+s_l+s_k))\Gamma(\frac12(1+s_l-s_k))
        \Gamma(\frac12(1-s_l+s_k))\Gamma(\frac12(1-s_l-s_k))}
 {\Gamma(\frac12(s_l+s_k))\Gamma(\frac12(s_l-s_k))
        \Gamma(\frac12(-s_l+s_k))\Gamma(\frac12(-s_l-s_k))}
          \times \\& \times
\Phi_{s_1,\dots,s_p}(t_1,\dots,t_p)\,ds_1ds_2\dots ds_p
\end{align}
where $C$  is a constant}.

\mal

{\scc Remark.} Factor (2.3)--(2.4) is the Gindikin--Karpelevich
density, it is an elementary function, see (1.44)--(1.46).

\mal

{\bff 2.\punkt. Analytic formula for arbitrary $\alpha$.}
Fix $m=0,1,\dots,p$. Consider nonnegative integers
  $$u_1\le u_2\le\dots\le u_m$$
satisfying the condition
  $$\alpha+2u_m+m<\frac12(p+q)$$
(if $m=0$, then a collection $\{u\}$ is empty).

{\scc Theorem 2.\fact.} {\it Let $p\ne q$ and
$\alpha$ be arbitrary, or $p=q$
 and $\alpha\in\R\setminus\{1,2,\dots,p-1\}$.
Then                         }
\begin{multline}
\prod_{k=1}^p\ch^{-\alpha} t_j=C\cdot
\sum\limits_{\begin{array}{c}m;\\ \quad u_1\le \dots\le u_m<
             \frac14(p+q)-\frac m2-\frac12\alpha\end{array}}
              E_m(\alpha, u)
\times\\ \times
\int\limits_{i\R^{p-m}} Y_m(\alpha,u; s)\r_m(s)
              \Phi_{\alpha-(p+q)/2+1+2u_1,\dots,
            \alpha-(p+q)/2+m+2u_m,s_{m+1},\dots,s_m}(t_1,\dots,t_p)
\times\\ \times
              ds_{m+1}\dots ds_p
\end{multline}
{\it where $C$ is the same as above},
\begin{align}
&E_m(\alpha,u)= (2\pi)^m\frac{p!}{m!}{2^\alpha} \prod_{j=1}^p\frac1{\Gamma(\alpha-j+1)}
                \times\\& \times
\prod\limits_{\tau=1}^m \frac
    {(-\alpha+\frac12(p+q)-2u_\tau-\tau)
     \Gamma(\alpha-p+\tau+2u_\tau)\Gamma(-\alpha+q-\tau-2u_\tau)}
    {(u_\tau-u_{\tau-1})!
      \Gamma(-\alpha+\frac12(p+q)-\tau+1-u_\tau-u_{\tau-1})}
\times\\& \times
\prod\limits_{1\le\sigma<\tau\le m} \Bigl\{
    (-\alpha+\half(p+q)-\half(\tau+\sigma)-u_\sigma-u_\tau)
     (\half(\tau-\sigma)+u_\tau-u_\sigma)
    \times\\& \times
    \frac{ \Gamma(\frac12(\tau-\sigma+1)+u_\tau-u_\sigma)
     \Gamma(-\alpha+\frac12(p+q)-\frac12(\tau+\sigma)
                   -u_\tau-u_{\sigma}+\frac12)
                           }
    {\Gamma(\frac12(\tau-\sigma)+u_\tau-u_{\sigma-1})
     \Gamma(-\alpha+\frac12(p+q)-\frac12(\sigma+\tau)-u_\sigma-u_\tau)}
\Bigr\}
\end{align}

\begin{align}&
Y_m(\alpha,u;s)=\notag
\\&
=\prod\limits_{k=m+1}^p\Bigl\{
    \Gamma(\half(\alpha-\half(p+q)+m+1+s_k))
 \Gamma(\half(\alpha-\half(p+q)+m+1-s_k))\Bigr\}
\times\\& \times
\prod\limits_{\tau\le m;\,\,k>m}
 \Biggl\{(\half(-\alpha+\half(p+q)-\tau-2w_\tau+s_k))
 (\half(-\alpha+\half(p+q)-\tau-2w_\tau-s_k))
\times\\& \times
\frac{\Gamma(\frac12(-\alpha+\frac12(p+q)-(\tau-1)-2u_\tau+s_k)
\Gamma(\frac12(-\alpha+\frac12(p+q)-(\tau-1)-2u_\tau-s_k)}
{\Gamma(\frac12(-\alpha+\frac12(p+q)-\tau+2+2w_{\tau-1}+s_k)
\Gamma(\frac12(-\alpha+\frac12(p+q)-\tau+2+2w_{\tau-1}-s_k)}\Biggr\}
\end{align}
{\it and  }
\begin{align}&
\r_m(s)=\prod_{k=m+1}^p\frac{\Gamma((q-p)/2+s_k) \Gamma((q-p)/2-s_k)}
                       {\Gamma(s_k)\Gamma(-s_k)}
\times\\& \times
\prod_{m+1\le k<l\le p}\frac
        {\Gamma(\frac12(1+s_l+s_k))\Gamma(\frac12(1+s_l-s_k))
        \Gamma(\frac12(1-s_l+s_k))\Gamma(\frac12(1-s_l-s_k))}
 {\Gamma(\frac12(s_l+s_k))\Gamma(\frac12(s_l-s_k))
        \Gamma(\frac12(-s_l+s_k))\Gamma(\frac12(-s_l-s_k))}
\end{align}

{\scc Remarks.} a) The factor $\r_m(s)$ is an elementary function.

   b) More convenient notations are used in  Section 4 (see 4.13).

   c) The formula, which is not
so explicit, but short is given in Section 6.

\mal

{\scc Remark.} The summand corresponding $m=0$ coincides with
integral (2.1)--(2.5). For summands corresponding $m=p$, the integration
is given by one point set and hence these summands are spherical
functions $\Phi_{\dots}$ with some coefficients.

\mal

{\bff 2.\punkt.  The case $\alpha=p-1, p-2,\dots, 1$.}
In this case some  summands disappear.

\mal

{\scc Proposition 2.\fact.} {\it Let $\alpha=p-h$ where $h\le p$.
Then the factor $E_m(\alpha, u)$ is nonzero if and only if
$$ m\ge h;\qquad u_1=u_2=\dots=u_h=0$$}

\mal

{\sc
 Proof.} Vanishing of $E_m(\alpha, u)$
is completely defined by a behavior of the factor
\begin{equation}
\frac{\prod_{\tau=1}^m \Gamma(\alpha-p+\tau+2u_\tau)}
     {\prod_{j=1}^p \Gamma(\alpha-j+1)}
\end{equation}
The denominator has a pole of  order $h$ at $\alpha=p-h$.
If the fraction is non-vanishing, then the numerator has a pole of
order $h$ at the same point.   \hfill$boxtimes$

\mal

{\bff 2.\punkt. The case $\alpha=-1,-2,-3,\dots$.}
Assume  $E_m(\alpha,u)\ne 0$. The denominator of
(2.16) has a pole of order $p$ in $\alpha$.
 Hence, the numerator also has a pole of order
$p$. Hence,
$$m=p$$
This means that all integrals in  Plancherel formula (2.6)
 vanish and we have only finite sum
of spherical functions with some coefficients.
The coefficient $E_m(\alpha,u)$ is nonzero iff
$$m+2u_m\le -\alpha$$

\mal

{\bff 2.\punkt. The Plancherel formula for the kernel representations
$T_\alpha$ of $\O(p,q)$.}

\mal

{\scc Theorem 2.\fact.}
  {\it Let $\alpha$ satisfies  Berezin conditions  {\rm(1.9)}.
Then}

 a){\it if $E_m(\alpha, u)\ne 0$ (see Subsection {\rm 2.3}), then
all spherical functions
 $$ \Phi_{\alpha - (p+q)/2+1+2u_1,\dots,
\alpha-(p+q)/2+m+u_m,s_{m+1},\dots,s_m}$$
  are positive definite.}

 b){\it  formula {\rm(2.6)-(2.15)} is really the Plancherel formula}

{\bff 2.\punkt.  The Plancherel formula for kernel-representations
of $\O(p+q)$.}
For a negative integer $\alpha=-n$
 (see Subsection 2.4 above)
formula (2.6) gives the expansion of $\det(1-z^*z)^n$ in
$\O(p)\times\O(q)$-spherical functions of $\O(p,q)$
and this is equivalent to the Plancherel formula for the kernel
representations of $\O(p+q)$.

\mal

{\bff 2.\punkt. The case of indefinite Shapovalov form.}
 For noninteger $\alpha<p-1$ we obtain
the problem discussed in Section 1.E.

\bol

\addtocounter{sec}{1}
\setcounter{equation}{0}
\setcounter{fact}{0}
\setcounter{punkt}{0}

{\bfff 3. B-function of the space
$\O(p,q)/\O(p)\times\O(q)$}

\nopagebreak

\bol

In this section we construct a matrix imitation of the $\B$-integral
$$\B(x,y)=\int_0^\infty\frac{t^{x-1}}{(1+t)^{x+y}}dt$$
for the symmetric spaces  $\O(p,q)/\O(p)\times\O(q)$.
For symmetric spaces $\GL(n,\K)/\U(n,\K)$ the $\B$-integrals
were defined by Gindikin \cite{Gin1} (see also exposition in \cite{FK}),
for other symmetric spaces $\B$-integrals were derived in \cite{Ner4}.

{\bff 3.\punkt. B-integral.}
Let
 $$\lambda_1,\dots,\lambda_p,\sigma_1,\dots,\sigma_p\in \C $$
 We also assume
$$\lambda_{p+1}=\sigma_{p+1}=0$$
Let $\SW_{p,q}$ be the section of wedge defined in Subsection 1.20.

{\scc Theorem 3.\fact.} {\it Let $\lambda_k$, $\sigma_k$ satisfy
the inequalities
\begin{equation}
\frac12(q+k)/2+1<\lambda_k<\sigma_k-\frac12(p-k)
\end{equation}
Then}
\begin{align}&
\int\limits_{\SW_{p,q}(\R)}
\prod_{j=1}^p
\frac{ \det[M-LL^t]_j^{\lambda_j-\lambda_{j+1}}  }
  {\det [1+M+N]_j^{\sigma_j-\sigma_{j+1}} }
  \cdot
\det(M-LL^*)^{-(p+q)/2} dM\,dN\,dL
=\\&=
\int\limits_%
{\begin{array}{c}M-LL^t>0\\ N=-N^t\end{array}}
\prod_{j=1}^n
\frac{
\det\left[\begin{array}{cc}1&L^t\\L&M\end{array}\right]_{q-p+j}%
^{\lambda_j-\lambda_{j+1}}    }
     {
\det[1+M+N]_j^{\sigma_j-\sigma_{j+1}}  }
  \det \left(\begin{array}{cc}1&L^t\\L&M\end{array}\right)^{-(p+q)/2}
 dL\,dM\,dN
=\\=&
=\prod_{k=1}^p
\pi^{k-(q-p)/2-1}
\frac{\Gamma(\lambda_k-(q+k)/2+1)\Gamma(\sigma_k-\lambda_k-(p-k)/2)}
     {\Gamma(\sigma_k-p+k)}
\end{align}

The proof of the Theorem is given in Subsections 3.2-3.6.

\mal

{\scc Remark.} For $p=q$ we have $L=0$ and  integral (3.2)--(3.3)
  has more
simple form, see (0.5). In this case the calculation given below also is simpler.
 The main simplification is the expression for  matrix (3.13):
the first block row and
the first block column are lacked.

\mal

{\scc Remark.}  We have $M=M^t>0$, $N=-N^t$. Hence,
$$\det(1+M+N)> 0$$
Indeed, for any $v\in\C^p$ we have $\Re v(M+N)v^*=vMv^*>0$. Hence, the
eigenvalues $\lambda_j$ of $M+N$ satisfy the condition
$\Re \lambda_j>0$. Hence, the eigenvalues of $1+M+N$ are nonzero.

\mal

{\scc Remark.} Hua Loo Keng in \cite{Hua}
evaluated the integrals%
\footnote{Hua integrals also can be reduced to the Selberg
$\B$-integrals by integration over $K\setminus G/K$.}
\begin{equation}
\int_{\B_{p,q}(\R)}\det(1-zz^*)^\tau dz
\end{equation}
Cayley transform reduces the {\it Hua integral} to the
following partial case of our integral
$$\const\cdot\int\limits_{SW_q(\R)}
\frac{\det\begin{pmatrix}1&L^t\\L&M\end{pmatrix}^{\tau}}
{\det(1+M+N)^{2\tau}}dL\,dM\,dN $$
Our calculation in this case is not homotopic to Hua calculations.

\mal

{\bff 3.\punkt. Replacement of notations.}
Firstly, we call to mind the standard formula (see \cite{Gan})
 for determinant of block
$(m+n)\times(m+n)$-matrix
\begin{equation}
\det\begin{pmatrix} A&B\\C&D\end{pmatrix}=
\det A \cdot \det (D-BA^{-1}C)
\end{equation}

Let us represent $M,N$ as  block
 $((p-1)+1) \times((p-1)+1 )$ matrices,  and $L$ as a block
 $((p-1)+1 ) \times (q-p)$ matrix:
$$
M=\left(\begin{array}{cc}P&q^t\\q&r\end{array}\right);\qquad
N=\left(\begin{array}{cc}A&-b^t\\b&0\end{array}\right)  ;\qquad
L=\left(\begin{array}{c} H\\l\end{array}\right)
$$
Then for $j\le p-1$
\begin{gather*}
\left[\begin{array}{cc}1&L^t\\L&M\end{array}\right]_{q-p+j}
\mbox{coincides with}
\left[\begin{array}{cc}1&H^t\\H&P\end{array}\right]_{q-p+j}
\\   \vphantom{\Biggl|^2}
[1+M+N]_j\qquad\mbox{coincides with}\qquad[1+P+A]_j
\end{gather*}
and by (3.6)
\begin{eqnarray}
\det \left(\begin{array}{cc}1&L^t\\L&M\end{array}\right)=
\det \left(\begin{array}{ccc}
     1&H^t&l^t \\
     H&P&q^t\\
     l&q&r  \end{array}\right) =
\qquad\qquad\qquad\qquad\qquad\qquad\nonumber\\
=\det \left(\begin{array}{cc}1&H^t\\H&P\end{array}\right) \cdot
\biggl[r -
 \left(\begin{array}{cc}l&q\end{array}\right)
    \left(\begin{array}{cc}1&H^t\\H&P\end{array}\right)^{-1}
     \left(\begin{array}{cc}l^t\\q^t\end{array}\right)
      \biggr]\nonumber \\
\det(1+M+N)=\det(1\!+\!P\!+\!A)
      \cdot\bigl(1+r-(q+b)(1\!+\!P\!+\!A)^{-1}(q^t-b^t)\bigr) \nonumber
\end{eqnarray}

By the Sylvester criterion the condition
$\begin{pmatrix}1&L^t\\L&M\end{pmatrix}>0$
(see (1.35)) in new notations has the form
\begin{equation}
\left(\begin{array}{cc}1&H^t\\H&P\end{array}\right)>0;\qquad
r -
 \left(\begin{array}{cc}l&q\end{array}\right)
    \left(\begin{array}{cc}1&H^t\\H&P\end{array}\right)^{-1}
     \left(\begin{array}{cc}l^t\\q^t\end{array}\right)>0
\end{equation}

By the remark given in Subsection 3.1
$$\det(1+P+A)>0$$

{\bff 3.\punkt. Substitution.}
Let us replace the variable $r$ to the variable
$$u=r -
 \left(\begin{array}{cc}l&q\end{array}\right)
    \left(\begin{array}{cc}1&H^t\\H&P\end{array}\right)^{-1}
     \left(\begin{array}{cc}l^t\\q^t\end{array}\right)
$$
(all other variables are the same). By (3.7) we have $u>0$.
 The Jacobian
of the substitution is 1.
Our integral converts to the form
\begin{align}
&\int dP\,dA\,dH\Biggl(\Xi (A,P,H)
 \times  \\
&\!\!\!\times\int\limits_{u>0,\,q,b\in\R^{p-1},l\in\R^{q-p}} \!\!\!\!\!\!\!\!\!\!\!
u^{\lambda_p-(p+q)/2}\biggl\{1\!+\!u\! +\!
 \left(\begin{array}{cc}l&q\end{array}\right)
    \left(\begin{array}{cc}1&H^t\\H&P\end{array}\right)^{-1}
     \left(\begin{array}{cc}l^t\\q^t\end{array}\right)+ \\
&\!\!\!\!+ \left(\begin{array}{cc}q&b\end{array}\right)
 \left(\begin{array}{cc}
    -(1\!+\!P\!+\!A)^{-1}& -(1\!+\!P\!+\!A)^{-1}\\
      (1\!+\!P\!+\!A)^{-1} & (1\!+\!P\!+\!A)^{-1} \end{array}\right)
     \left(\begin{array}{cc}q^t\\b^t\end{array}\right)
\biggr\}^{-\sigma_p}
du\,dl\,dq\,db\Biggr)
\end{align}
where
\begin{equation}
\Xi (A,P,H)=\prod_{j=1}^{p-2}
\frac{
\det\left[\begin{array}{cc}1&H^t\\H&P\end{array}\right]_{q-p+j}%
^{\lambda_j-\lambda_{j+1}}    }
     {
\det[1+P+A]_j^{\sigma_j-\sigma_{j+1}}  } \cdot
 \frac{ \det \left(\begin{array}{cc}1&H^t\\H&P\end{array}\right)%
^{\lambda_{n-1}-(p+q)/2}}
{\det(1+P+A)^{\sigma_{n-1}}  }
\end{equation}
is an expression independent on $u,b,l,q$.

Firstly, we want to evaluate interior integral (3.9)--(3.10)

\mal

{\bff 3.\punkt.
 Transformation of the integrand.} Denote by $S$ the expression
$$S=1+P+A$$
Let us represent
the expression in the curly brackets in (3.10)  in the form
\begin{equation}
\biggl\{1+u+
 \left(\begin{array}{ccc}l&q&b\end{array}\right)
                            X
 \left(\begin{array}{c}l^t\\q^t\\b^t\end{array}\right)
\biggr\}
\end{equation}
where
\begin{equation}
X=
     \left(\begin{array}{ccc}
    \left(\begin{array}{cc}1&H^t\\H&P\end{array}\right)^{-1}+
    \left(\begin{array}{cc}0&0\\0&-S^{-1}\end{array}\right)
                        &   \phantom{a}&
 \left(\begin{array}{c}0\\-S^{-1}\end{array}\right)
                \\
  \phantom{d} & & \\
     \left(\begin{array}{cc}0& S^{-1}\end{array}\right)
             &   &
            S^{-1}
\end{array}\right)
\end{equation}
(we wright a block matrix whose elements are block
matrices itself).
The last summand in the curly brackets is
a quadratic form in the variables $b,q,l$. But the matrix $X$
is not symmetric and it is more natural to
re-wright  expression (3.12) in the form
\begin{equation}
\biggl\{1+u+
 \left(\begin{array}{ccc}l&q&b\end{array}\right)
                            \frac12(X+X^t)
 \left(\begin{array}{c}l^t\\q^t\\b^t\end{array}\right)
\biggr\}
\end{equation}

{\bff 3.\punkt. Separation of variables.}

{\scc Lemma 3.\fact.}
$$\det\left(\frac12(X+X^t)\right)=
\det    \left(\begin{array}{cc}1&H^t\\H&P\end{array}\right)^{-1}
\cdot
\det (1+P+A)^{-2}$$

{\scc Proof.}
$\det\left({\textstyle\frac12}(X+X^t)\right) =$
$$=\det
     \left(\begin{array}{cc}
    \left(\begin{array}{cc}1&H^t\\H&P\end{array}\right)^{-1}+
    \left(\begin{array}{cc}0&0\\0&-\frac12S^{t-1}-\frac12S^{-1}\end{array}\right)
                        &
 \left(\begin{array}{c}0\\-\frac12S^{-1}+\frac12S^{t-1}\end{array}\right)
                \\
  \phantom{d} &  \\
     \left(\begin{array}{cc}0& \frac12S^{-1}-\frac12S^{t-1}\end{array}\right)
             &
           \frac12 S^{-1}+\frac12S^{t-1}
\end{array}\right)
$$
Adding the third row to the second row
and the third column to the second column,
we obtain
$$\det
     \left(\begin{array}{cc}
    \left(\begin{array}{cc}1&H^t\\H&P\end{array}\right)^{-1}
                        &
 \left(\begin{array}{c}0\\(1+P-A)^{-1}\end{array}\right)
                \\
  \phantom{d} &  \\
     \left(\begin{array}{cc}0& (1+P+A)^{-1}\end{array}\right)
             &
           \frac12 (1+P+A)^{-1} +\frac12 (1+P-A)^{-1}
\end{array}\right)
$$
Formula (3.6) reduces the determinant to the form
\begin{eqnarray*}\det
    \left(\begin{array}{cc}1&H^t\\H&P\end{array}\right)^{-1}
      \cdot
\det\Bigl( \frac12 (1+P+A)^{-1}  + \frac12 (1+P-A)^{-1}  - \\-
 \left(\begin{array}{cc}0& (1+P+A)^{-1}\end{array}\right)
  \left(\begin{array}{cc}1&H^t\\H&P\end{array}\right)
 \left(\begin{array}{c}0\\(1+P-A)^{-1}\end{array}\right)
    \Bigr)
=    \\  =
  \det  \left(\begin{array}{cc}1&H^t\\H&P\end{array}\right)^{-1}
\det(1+P+A)^{-1}    \det (1+P-A)^{-1}   \times\\ \times
\det\bigl({\textstyle\frac12} (1+P-A) + {\textstyle\frac12}(1+P+A) -P\bigr)
\end{eqnarray*}

The last factor is 1. We also observe
$$(1+P+A)^t=1+P-A$$
and hence their determinants coincides.   \hfill $\boxtimes$

\mal

{\scc Lemma 3.\fact.} $X+X^t>0$.

\mal

{\scc Proof.}  In the identity
$$
\det(1+M+N)=\det(1\!+\!P\!+\!A)
      \cdot\bigl[1+r-(q+b)(1\!+\!P\!+\!A)^{-1}(q^t-b^t)\bigr] \nonumber
$$
we have $\det(1+M+N)>0$, $\det(1+P+A)>0$. Hence, the factor in the square
brackets is positive. Hence, expression (3.12) is positive for all
$u>0$, and all  $q,b,l$. Quantity (3.12) coincides with quantity (3.14).
Hence, the matrix $X+X^t$ is nonnegative defined. By Lemma 3.2
its determinant is nonzero and we obtain the required statement.
\hfill $\boxtimes$

\mal

Consider  the linear substitution
$$
 \left(\begin{array}{ccc}l&q&b\end{array}\right)
\sqrt{{\textstyle\frac12}(X+X^t)}=h\in\R^{q-p}\oplus\R^{p-1} \oplus\R^{p-1}
$$
to  interior integral (3.9)--(3.10). Its Jacobian is
$$\det\begin{pmatrix}1&H^t\\H&P\end{pmatrix}^{1/2}\cdot
\det(1+P+A)
$$
and hence the interior integral coverts to the form
\begin{align}&
\det\begin{pmatrix}1&H^t\\H&P\end{pmatrix}^{1/2}\cdot
\det(1+P+A)\times \\  \times
&
\int\limits_{u>0,\,h\in\R^{q+p-2}}
u^{\lambda_p-(p+q)/2}\bigl\{1+u +  |h|^2
\bigr\}^{-\sigma_p}
du\,dh
\end{align}
The first factor (3.15) adds to the product $\Xi(A,P,H)$ (see (3.11))
and we reduce our $\B$-integral   (3.3) to the product of the integrals
\begin{multline*}
\int_{P-HtH>0,\,A=-A^t} \Xi(A,P,H)
\det\begin{pmatrix}1&H^t\\H&P\end{pmatrix}^{1/2}\cdot
\det(1+P+A)dA\,dP\,dH \times\\ \times
\int\limits_{u>0,\,h\in\R^{q+p-2}}
u^{\lambda_p-(p+q)/2}\bigl\{1+u +  |h|^2
\bigr\}^{-\sigma_p}
du\,dh
\end{multline*}

Let us denote  $\B$-integral (3.3)
by
\begin{equation}
I_{p,q}(\alpha_1,\dots,\alpha_p;\sigma_1,\dots,\sigma_p)
\end{equation}
and let us denote factor (3.16) by
$J_{p,q}(\alpha_p;\sigma_p)
$.
We obtain the following recurrent identity
\begin{multline*}
I_{p,q}(\alpha_1,\dots,\alpha_p;\sigma_1,\dots,\sigma_p)=
\\=
I_{p-1,q-1}(\alpha_1,\dots,\alpha_{p-2},\alpha_{p-1}-\half;
\sigma_1,\dots,\sigma_{p-2},\sigma_{p-1}-1)
J_{p,q}(\alpha_p;\sigma_p)
\end{multline*}

{\bff 3.\punkt. Evaluation of $ J_{p,q}(\alpha_p;\sigma_p)$.}
This problem is trivial. Firstly, we consider spherical coordinates
in $\R^{p+q-2}$ in the variable $h$. Then
$ J_{p,q}(\alpha_p;\sigma_p)$ converts to the
form
$$\frac{2\pi^{(p+q)/2-1}}{\Gamma((p+q)/2-1)}
\int_{u>0}\int_{r>0}u^{\lambda_p-(p+q)/2}r^{p+q-3}\bigl\{1+u +  r^2
\bigr\}^{-\sigma_p}
dr\,du
$$
The substitution $v=r^2$ reduces our integral to a special case of
the Dirichlet $\B$-integral
 $$\int_{u>0, v>0}
\frac{ u^{a-1}v^{b-1}}
     {(1+u+v)^{a+b+c}}du\,dv=
\frac{\Gamma(a)\Gamma(b)\Gamma(c)}{\Gamma(a+b+c)}$$

This completes the proof of Theorem 3.1.

\mal

{\bff 3.\punkt. Spherical transform of $\b_\alpha$.}

\mal

{\scc  Corollary 3.\fact.} {\it Let $\alpha>p+q-1$.
Then spherical transform of $\b_\alpha$ is}
\begin{equation}
\frac{2^\alpha}{\prod_{1\le j\le p}\Gamma(\alpha-j+1)}
\prod\limits_{k=1}^p\Gamma(\half(\alpha-\half(p+q)+1+s_k))
                    \Gamma(\half(\alpha-\half(p+q)+1-s_k))
\end{equation}

{\scc Proof.} The function $\b_\alpha$ is given by the formula (1.37).
By Subsection 1.24 we must evaluate the integral
$$\int_{G/K}\b_\alpha(z)\Psi_s(z)\,d\lambda(z)$$
But the integral is a special case of our $\B$-integral.

        \mal

{\bff 3.\punkt. Proof of Theorem 2.1.}
By the Gindikin-Karpelevich inversion formula and Corollary 3.4
we obtain the statement of the theorem for $\alpha>(p+q)-1$.

For   $\alpha>(p+q)/2-1$ the statement of the Theorem follows
from trivial Lemma 4.1 proved below.

\bol

{\bfff 4. Formal analytic continuation.}

\nopagebreak

\bol

\addtocounter{sec}{1}
\setcounter{equation}{0}
\setcounter{fact}{0}
\setcounter{punkt}{0}

We proved the Plancherel
 formula (2.1)-(2.5) for large values of the parameter $\alpha$.
Its left part $\prod\ch^{-\alpha}(t_j)$
 depends analytically on $\alpha\in\C$.
The integrand in the right part has
 singularities on the lines
\begin{equation}
\Re\alpha=\half(p+q)-1-2\kappa; \qquad \mbox{where}\quad
 \kappa=0,1,2,\dots
\end{equation}
Thus, the right part of formula (2.1)--(2.5)
may be nonalytic for these values of $\alpha$.

 Our
next purpose is to construct the analytic continuation of
the right part
to arbitrary complex $\alpha$.

\mal

{\bff 4.\punkt. Analyticity.} Let us denote the right part of the
 formula (2.1)--(2.5)
by
\begin{equation}
{\frak F}(\alpha):={\frak F}(\alpha;t)=E(\alpha)
\int_{i\R^p} Y(\alpha;s)\r(s)\Phi_s(t)\,ds
\end{equation}
where the meromorphic factor $E(\alpha)$ is given by  formula (2.1),
the factor   $Y(\alpha;s)$ is defined by (2.2) and $\r(s)$
is  Gindikin--Karpelevich density (2.3)--(2.4).
{\it In this Section we  fix the variable $t$ and we  omit
the argument $t$ from  the notation   ${\frak F}(\alpha;t)$.}

Consider domains $\Pi_0$, $\Pi_1$, \dots in $\C$ defined by
\begin{align*}
&\Pi_0: \qquad \Re\alpha>\half(p+q)-1\\
&\Pi_k: \qquad  \half(p+q)-1-2k<
  \Re\alpha< \half(p+q)-1-2(k-1)\quad\mbox{where}\,\, k>0
\end{align*}

{\scc Lemma 4.\fact.} {\it The function ${\frak F}(\alpha)$
is an analytical function
on $\Pi_\kappa$ for all $\kappa=0,1,2\dots$.}

\mal

{\sc Proof.} a) {\it Convergence of  integral} (4.2).
First, the Gindikin-Karpelevich factor $\r(s)$ has a polynomial growth
in $s$, see formulas (1.44)--(1.46).

By the formula (see \cite{HTF},1.18.6)
\begin{equation}
|\Gamma(a+iy)|=(2\pi)^{1/2}|y|^{a-1/2}\exp\{\half \pi|y|\}(1+o(1));
 \qquad |y|\to\infty
\end{equation}
the factor $Y(\alpha;s)$ exponentially decreases.

A spherical function $\Phi_s(t)$ is a spherical function
of an unitary representation and hence we have $|\Phi_s(t)|\le 1$%
\footnote{For following inductive steps this arguments
must  be replaced by inequality (1.39)}.

Hence, the integrand exponentially decreases and the integral
absolutely converges.

\mal

b){\it Existence of $\frac \partial{\partial\alpha}{\frak F}(\alpha)$}.
It is sufficient to prove
uniform convergence of the integral
\begin{equation}
\int_{i\R^p}\frac\partial{\partial\alpha}Y(\alpha; s)\r(s)\Phi_s(t)ds
\end{equation}in small neighborhood of a fixed point
$\widetilde\alpha$. For this we needs in
uniformity by $a$ of $o(1)$ in (4.3). In fact, the
asymptotics is really uniform
but formally we have no possibility to refer to \cite{HTF}.
Formula
(4.3) is derived from the Binet formula (see \cite{HTF},(1.9.4))
\begin{multline*}
  \ln\Gamma(z)=(z-\frac12)\ln z-z+\frac12\ln(2\pi)
 + \int_0^\infty \left[\frac1{e^t-1}-\frac 1t+\frac12\right]t^{-1}e^{-tz}dt
\end{multline*}
This formula easily implies an uniform estimate of the form
$$|\Gamma(a+iy)|\le \const\cdot\exp(-(\half\pi-\epsilon)|y|);
          \qquad |a-\widetilde a|<\delta$$
The Cauchy integral for derivative
$$-2\pi i f'(z)=\int_L\frac {f(z)}{(z-u)^2}dz$$
implies the same estimate for derivative of $\Gamma$-function.

We observe that integrand (4.4) is dominated by some function
having the form
$$P(s)\exp(-b\sum |s_j|)$$
where $P(s)$ is a polynomial and $b>0$.

  Thus, the function ${\frak F}(\alpha)$ has a derivative in the
 complex variable
$\alpha$ and this complets the proof.  \hfill$\boxtimes$

\mal

{\scc Lemma 4.\fact} {\it Let $p\ne q$. Then the function $\F(\alpha)$
is continuous on the line $\alpha\in \R$.}

\mal

{\scc Proof.} Let $h=\frac12(p+q)-1-2\kappa$ be one of our singular points.
The singularity of the integrand near this point has the form
$$\const\cdot\prod_{k=1}^p \frac {s_k^2}{(\alpha-h)^2-s_k^2}
\prod_{1\le k<l\le p}(s_k^2-s_l^2)^2 (1+o(1)) $$
Remind that $s_k$ are pure imaginary.
Hence, the integrand is bounded in a neighborhood of the
 point $\alpha=h, s=0$.
As we have seen in the previous proof, the integrand has
an integrable majorant in a domain $|s_k|>A$, $|\alpha-h|<\epsilon$
($\alpha\in\R$). By Lebesgue theorem about dominant convergence,
 expression (4.2) is continuous at the point
$\alpha=h$.   \hfill$\boxtimes$

\mal

{\scc Remark.} The function $\F(\alpha)$ is continuous at the real points
$\alpha=\half(p+q)-1-2\kappa$ but it is not smooth at these points.

\mal

  We denote the restriction of the function $\F(\alpha)$ to the domain
$\Pi_\kappa$ by
$$\F_\kappa(\alpha)$$

\mal

{\bff 4.\punkt. Analytic continuation of $\F_\kappa(\alpha)$
 through a point of a line $\Re\alpha=\frac12(p+q)-1-2\kappa$.}
The following lemma is  main in this Section.
Its proof is given in Subsections 4.2-4.5.

\mal

{\scc Lemma 4.\fact.} {\it Let $\alpha_0$ satisfies the condition
$$\Re \alpha_0=\half(p+q)-1-2\kappa; \qquad \Im \alpha_0\ne 0$$
Then}

\mal

a) {\it the function $\F_\kappa(\alpha)$ admits the analytic continuation
to some small neighborhood
 \begin{equation}
{\cal O}_\delta: \qquad  |\alpha-\alpha_0|<\delta \qquad\mbox{where}\quad
\delta<\min\{1/1000, |\Im\alpha|/1000\}
\end{equation}
of the point $\alpha_0$ }

\mal

b) {\it for any
$$\alpha\in {\cal O}\cap \Pi_{\kappa+1}$$
we have}
\begin{align}
&\F_{\kappa+1}(\alpha)-\F_{\kappa}(\alpha)
=\frac{\pi p}{\kappa!}\prod_{j=1}^p \frac1{\Gamma(\alpha-j+1)}
\times \\&\times
\frac{(-\alpha+\frac12(p+q)-2\kappa-1)
      \Gamma(\alpha-p+1+2\kappa)
      \Gamma(-\alpha+q-1-2\kappa)
      }
     {\Gamma(-\alpha+\frac12(p+q)-\kappa)}
\times\\& \times
\int\limits_{i\R^{p-1}}
\prod\limits_{2\le k\le p, \pm}
  \Gamma(\frac12(\alpha-\frac12(p+q)+2+2\kappa\pm s_k))
\times\\& \times
 \prod\limits_{2\le k\le p,\pm}
\frac{ (\frac12(-\alpha+\frac12(p+q)-1\pm s_k))
       \Gamma(\frac12(-\alpha+\frac12(p+q)-2\kappa\pm s_k))}
        {\Gamma(\frac12(-\alpha+\frac12(p+q)+1\pm s_k))}
\times \\& \times
\prod\limits_{2\le k\le p, \pm}\frac {\Gamma(\frac12(q-p)\pm s_k)}
                          {\Gamma(\pm s_k)}
\times \\ & \times
\prod\limits_{2\le k<l\le p,\pm}
 \frac{\Gamma(\frac12(1+s_k\pm s_l))\Gamma(\frac12(1-s_k\pm s_l)}
{\Gamma(\frac12(s_k\pm s_l)\Gamma(\frac12(-s_k\pm s_l)}
\times\\ &\times
\Phi_{\alpha-(p+q)/2+1+2\kappa, s_2, \dots, s_p}(t)
     ds_2\dots  ds_p
\end{align}

In the last formula we use the following  notation
\begin{equation}
\prod\limits_{\pm}\Gamma(a\pm s):=\Gamma(a+s)\Gamma(a-s)
\end{equation}

\begin{figure}
\risodin

{\it Picture} 1. The complex plane $\alpha$. The lines
 $\alpha=(p+q)/2-1-2\kappa$ and the  domains $\Pi_k$
\end{figure}

    \mal

{\bff 4.\punkt. Existence of the analytic continuations.}
Let us represent  expression (4.2) (or (2.1)--(2.5))
for $\F_\kappa(\alpha)$ in the form
$$E(\alpha)
\int \prod_{k,\pm}\Gamma(\half(\alpha-\half(p+q)+1\pm s_k))\cdot
      \r(s)\Phi_s(t)ds$$

Let $\epsilon_1>\epsilon_2\dots>\dots>\epsilon_p\ge 0$
be very small (for instance $\epsilon_1<\delta/10$
where $\delta$ was defined in (4.5)).
Consider the function

\begin{multline}
\F_\kappa(\alpha;\epsilon)=
E(\alpha)\int_{i\R^p}
 \prod_{k,\pm}\Gamma(\frac12(\alpha+\epsilon_k-\half(p+q)+1\pm s_k))\cdot
      \r(s)\Phi_s(t)ds
\end{multline}
in the domain
$$\Pi_\kappa^\epsilon:\qquad
          -2\kappa < \Re(\alpha-\frac12(p+q)+1)<-2(\kappa-1)-\epsilon_1$$

{\scc  Lemma 4.\fact.} a) {\it The function $\F_\kappa(\alpha;\epsilon)$
admits the holomorphic continuation to the domain $\cal O_\delta$}
\rm(\it see \rm (4.5)).

\mal

 b){\it The functions $|\F_\kappa(\alpha;\epsilon)|$
in $\cal O_\delta$ are bounded by a constant independent on $\epsilon$.}

\mal

 {\sc Proof.}
a) The factor $\r(s)$ is holomorphic in the domain
$|\Re(s_j)|<1/4$ and its poles are very far from the contour $\L$
which is described below.

\begin{figure}
\risdva

Picture 2.
 The contour $L_k=R_1\cup S_1\cup Q\cup S_2\cup R_2$
 on the complex plane $s_k$. The centers
of the semicircles $S_1,S_2$ are $\pm\Im\alpha_0$, the radius of
the semicircles
is $10\delta$.
\end{figure}

 Consider
\begin{equation}
\alpha\in{\cal O}_\delta\cap\Pi_\kappa^\epsilon
\end{equation}
Then integrand (4.14) has poles on hyperplanes
$$s_k=\pm(\half(p+q)-1-\alpha-2u-\epsilon_k);
\qquad\mbox{where}\quad u=0,1,\dots$$
If $u=\kappa$, then the poles are lying near points $\pm\Im\alpha_0$,
on  Picture 2 the poles are marked as black circles. The arrows show
the direction of their motion if $\Re\alpha$ decreases.
The white circles show rough position of the poles than
$\alpha\in{\cal O}_\delta\cap\Pi_{\kappa+1}^\epsilon$

  Consider the contour $L_k$ on the complex plane $s_k\in \C$
given by Picture 2.
Let ${\cal L}\subset \C^p$ be the product
of the contours $L_k$. Obviously, for
$\alpha\in{\cal O}_\delta $,  we can replace the
integration over  $i\R^p$
in the formula (4.14)
by the integration over $\cal L$.
 But the integral
$$\int_{\cal L}=\int_{\cal L}
\prod_{k,\pm}\Gamma(\frac12(\alpha+\epsilon_k-\half(p+q)+1\pm s_k))\cdot
      \r(s)\Phi_s(t)ds
$$
 obviously is
holomorphic with respect to $\alpha$
 in the domain ${\cal O}_\delta$
(indeed, the surface $\cal L$ doesn't intersect the singularities
and the integrand exponentially decreases as $|s|\to\infty$).

\mal

b) Consider the parameter $\theta_k:=\Im s_k$ on the contour
 $L_k$.

\mal

{\scc Lemma 4.\fact.} {\it There exist constants $A=A(t)$, $N$ such that
$$| \prod_{k,\pm}\Gamma(\half(\alpha+\epsilon_k-\half(p+q)+1\pm s_k))\cdot
      \r(s)\Phi_s(t)|\le
  A\prod_{j=1}^p \dot(1+|\theta_j|)^N\exp(-\pi|\theta_j|)$$
for all $(s_1,\dots, s_p)\in {\cal L}$.}

\mal

{\scc Proof of Lemma 4.5.} Let us estimate all factors in the left part
of the inequality.

a) {\it The Gindikin--Karpelevich factor} $\r(s)$.
By the formulas (1.45)--(1.46) for any imaginary $s_k$ we have
$$
\Bigl|\frac{\Gamma(\frac12(q-p)+s_k)\Gamma(\frac12(q-p)-s_k)}
{\Gamma(s_k)\Gamma(-s_k)}\Bigr|\le \const\cdot
         (1+|\theta_k|)^{2(q-p)}
$$
The same expression is bounded on the semi-circles $S_1, S_2$.

We also must estimate the factor (1.44).
Firstly,
$$|\prod_{k=1}^p(s_k^2-s_l^2)|\le\const \prod_k(1+|\theta_k|)^{2(p-1)}$$

Secondly, let us estimate the factors
$$\tg(\pi(s_k\pm s_l))$$
of (1.44).
If $s_k,s_l$ are imaginary, then  $|\tg(\pi(s_k\pm s_l))|<1$.
If $s_k,s_l\in S_1, S_2$, then this expression is bounded
(since $S_1,S_2$ are compact sets).
Let $s_k$ be imaginary and $s_l\in S_1, S_2$. Then
we obtain a value having the form $|\tg(x+iy)|$
where $x,y\in\R$, $|x|<10\pi\delta$. Then
$$|\tg(x+iy)|=\bigl|\frac{\tg x +\tg iy}{1+\tg x\tg iy}\bigr|\le
|\tg x +\tg iy|\le |\tg x|+1\le \tg (10\pi\delta)+1$$

b) {\it The  $\Gamma$-factor $Y(\alpha;s)$.}
 By  formula (4.3) for imaginary $s_k$
we have
\begin{multline*}
|\Gamma(\half(\alpha+\epsilon_k-\half(p+q)+1+ s_k))
     \Gamma(\half(\alpha+\epsilon_k-\half(p+q)+1- s_k)|\le
\\ \le
   \const\cdot (1+|\theta_k|)^{\Re\alpha+\epsilon_k-\frac12(p+q)+1}
               \exp(-\pi|\theta_k|)
\end{multline*}
For $s_k\in S_1, S_2$ the same expression is  bounded
(but very large).

c){\it Spherical functions $\Phi_s(t)$.}
By estimation (1.39) we have
$$|\Phi_s(t)|\le\Phi_{\Re s}(t)\le
\max \Phi_{r_1,\dots,r_p}(t)
$$
where maximum is given over all real vectors $(r_1,\dots,r_p)$
satisfying the condition $|r_j|\le10\delta$.
Hence, for a fixed $t$ the spherical function in integrand
is dominated by a constant.

This  completes the proof of the Lemma 4.5.\hfill $\boxtimes$

Now we can complete the proof of Lemma 4.4.b).
By Lemma 4.5 we have
$$|\F_\kappa(\alpha;\epsilon)|\le
  A(t) \int_{\R^p}\prod_{j=1}^p \bigl[
(1+|\theta_j|)^N\exp(-\pi|\theta_j|)\bigr]
\prod_{j=1}^p \chi(\theta_j) d\theta_1\dots d\theta_p  $$
there the function $\chi(\theta)$
is given by the formula
$$\chi(\theta):=\frac{ds}{d\theta}=\Bigl\{
 \begin{array}{ll} 1 & \quad \mbox{if}\quad
      \bigl|\Im\alpha_0-|\theta|\bigr|\ge 10\delta\\
  (100\delta^2-(\Im\alpha_0-\theta)^2)^{-1/2}
&\quad \mbox{if}\quad
      \bigl|\Im\alpha_0-|\theta|\bigr|\le 10\delta
\end{array}
$$
Hence
$$|\F_\kappa(\alpha;\epsilon)|\le A(t)^p \left(\int_{-\infty}^{\infty}
(1+|\theta|)^N\exp(-\pi|\theta|)\chi(\theta)d\theta\right)^p
$$
This completes the proof of uniform boundedness of the functions
$\F_\kappa(\alpha;\epsilon)$ for a fixed $t$
 (Lemma 4.4b). \hfill $\boxtimes$

\mal

Now we are ready to prove existence of the analytic continuation
of the function $\F_\kappa$.

\mal

{\sc Proof of Lemma 4.3.\rm b}. Let us denote by $\epsilon/n$ the vector
$(\epsilon_1/n,\dots,\epsilon_p/n)$.
Consider the sequence of functions
$$g_n(\alpha)=\F_\kappa(\alpha;\epsilon/n)$$
in the circle ${\cal O}_\delta$.
Since the functions $g_n(\alpha)$ are uniformly bounded,
by Montel theorem
there exists a subsequence $g_{n_j}$ which is uniformly
convergent on each smaller circle.
Let $g(\alpha)$ be its limit.
By Weierstrass theorem $g(\alpha)$ is holomorphic in ${\cal O}_\delta$.
It remains to notice that
 $$\lim_{n\to\infty}\F_\kappa(\alpha;\epsilon/n)=\F_\kappa(\alpha)$$
for $\alpha\in \Pi_\kappa\cap{\cal O}_\delta$.
Hence, $g(\alpha)$ is the  analytic continuation of $\F_\kappa(\alpha)$
to the circle ${\cal O}_\delta$.

\mal

{\bff 4.\punkt. Forcing of poles.} First, we want to obtain
an explicit formula for
the analytic continuation of $\F_\kappa(\alpha;\epsilon)$
to the domain $\Pi_{\kappa+1}^\epsilon$.

Let the contours $L_k$ be the same as above. Let $i\R_k$ be the imaginary axis
on the complex plane $s_k$. Consider the surface
$${\cal L}_k=i\R_1\times\dots\times i\R_{k-1}\times L_k\times\dots\times L_p
\subset\C^p$$
We have ${\cal L}_1=\cal L$, ${\cal L}_p=i\R^p$.

Consider $\alpha\in{\cal O}_\delta\cap \Pi_{\kappa+1}$.
Then
\begin{align}
&\F_{\kappa+1}(\alpha;\epsilon)=
E(\alpha)\int\limits_{i\R^p}
 \prod_{k,\pm}\Gamma(\half(\alpha+\epsilon_k-\half(p+q)+1\pm s_k))\cdot
      \r(s)\Phi_s(t)ds\\
&\F_\kappa(\alpha;\epsilon)=
E(\alpha)\int\limits_{\cal L}
 \prod_{k,\pm}\Gamma(\half(\alpha+\epsilon_k-\half(p+q)+1\pm s_k))\cdot
      \r(s)\Phi_s(t)ds
\end{align}
Hence,
\begin{equation}
\F_{\kappa+1}(\alpha;\epsilon)-\F_{\kappa}(\alpha;\epsilon)=
\int\limits_{i\R^p}-\int\limits_{\cal L}=
\sum_\sigma \Bigl[\int\limits_{{{\cal L}_{\sigma+1}}}-
       \int\limits_{{\cal L}_\sigma}\Bigr]
\end{equation}
Looking to Picture 2 we observe
\begin{align}
&\int_{{{\cal L}_{\sigma+1}}}-\int_{{\cal L}_\sigma}=\\
&=2\pi i
\int\limits_%
  {\begin{array}{c}s_1\in i\R_1,\dots,s_{\sigma-1}\in i\R_{\sigma-1}\\
s_{\sigma+1}\in L_{\sigma+1},\dots,\sigma_p\in L_p\end{array}}
\Bigl[\Res\limits_{s_\sigma=\alpha+\epsilon_\sigma-\frac12(p+q)+1+2\kappa}-
\\&-\Res\limits_{s_\sigma=
-\alpha-\epsilon_\sigma+\frac12(p+q)-1-2\kappa}\Bigr]
ds_1\dots ds_{\sigma-1}ds_{\sigma+1}\dots ds_p
\end{align}
The integrand in (4.16)--(4.17) is an
 even function in $s_\sigma$ and hence two residues
in (4.20)--(4.21)
differ only by  sign.
The order of poles
$$s_\sigma=\pm(\alpha+\epsilon_\sigma-\half(p+q)+1+2\kappa)$$
 of the integrand is 1 and hence the residues can be evaluated by
a simple substitution
\begin{multline}
\Res\limits_{s_\sigma=
\alpha+\epsilon_\sigma-\frac12(p+q)+1+2\kappa}=
      H_\sigma(\alpha,\epsilon,s):=\\=
E(\alpha)\biggl[
\frac{\Gamma(\frac12(\alpha+\epsilon_\sigma-\frac12(p+q)+1-s_\sigma))}
{s_\sigma-\alpha-\epsilon_\sigma+\frac12(p+q)-1-2\kappa}
\Gamma(\half(\alpha+\epsilon_\sigma-\half(p+q)+1+s_\sigma))
\times\\ \times
\prod_{\pm; k\ne\sigma}
    \Gamma(\half(\alpha+\epsilon_k-\half(p+q)+1\pm s_k))\cdot
      \r(s)\Phi_s(t)\biggr]
\Biggr|_{s_\sigma=\alpha+\epsilon_\sigma-\frac12(p+q)+1+2\kappa}
\end{multline}
In this way, we reduce sum (4.18) to
$$2\sum\limits_{\sigma=1}^p
\int\limits_{i\R_1\times\dots\times
    i\R_{\sigma-1}\times L_{\sigma+1}\times\dots\times L_p}
H_\sigma(\alpha,\epsilon,s)ds_1\dots ds_{\sigma-1}ds_{\sigma+1}\dots ds_p$$

We obtain an expression for
$\F_{\kappa+1}(\alpha;\epsilon)-\F_{\kappa}(\alpha;\epsilon)$.
Unfortunately the domains of integration is yet
complicated. By this reason we apply transformation (4.18) to each
summand in the last expression. We obtain $p(p-1)/2$ additional summands
which are
integrals over $(p-2)$-dimensional surfaces. Each integral can be easily
evaluated by residues. After this  we
 apply our arguments again, again, again.

It is possible to wright the final expression (it is slightly long).
 Fortunately,
this is not necessary.
The only goal of our interest is
\begin{equation}
\lim_{\epsilon\to 0}
(\F_{\kappa+1}(\alpha;\epsilon)-\F_{\kappa}(\alpha;\epsilon))
\end{equation}
For instance, consider  the summand obtained by the substitution
\begin{align*}
&s_\sigma=\alpha+\epsilon_\sigma-\frac12(p+q)+1+2\kappa\\
&s_\upsilon=\alpha+\epsilon_\upsilon-\frac12(p+q)+1+2\kappa
\end{align*}
Then the integrand contains the  factors
$$\frac1{\Gamma(\pm\frac12(s_\sigma-s_\upsilon))}\Bigr|_%
 {\begin{array}{l}
s_\sigma=\alpha+\epsilon_\sigma-\frac12(p+q)+1+2\kappa\\
s_\upsilon=\alpha+\epsilon_\upsilon-\frac12(p+q)+1+2\kappa
\end{array}}=
\frac1{\Gamma(\pm\frac12(\epsilon_\sigma-\epsilon_\upsilon))}
$$
This factors tend to 0 if $\epsilon\to 0$.
Of course, it is necessary to check lacking of poles of the numerator in
dangerous for us domain.

 Therefore a nonzero contribution to the limit
(4.23) can be given only by the terms
\begin{equation}
2\sum\limits_{\sigma=1}^p
\int\limits_{i\R^p}
H_\sigma(\alpha,\epsilon,s)ds_1\dots ds_{\sigma-1}ds_{\sigma+1}\dots ds_p
\end{equation}

Our expression is symmetric under permutations of $s_j$
and hence all summands of (4.28) give the same contribution to the limit.
Thus, deleting $\epsilon$ in (4.24) we obtain the formula
\begin{align*}
&\F_{\kappa+1}(\alpha)-\F_{\kappa}(\alpha)=\\
&\frac{2\pi i\cdot 2p\cdot 2^\alpha}{\prod\Gamma(\alpha-j+1)}
\int_{\R^{p-1}}
\biggl[\frac{\Gamma(\frac12(\alpha-\half(p+q)+1-s_1))}
{s_1-\alpha+\frac12(p+q)-1-2\kappa}
\Gamma(\half(\alpha-\half(p+q)+1+s_1))
\times\\& \times
\prod_{2\le k\le p,\pm}\Gamma(\half(\alpha-\half(p+q)+1\pm s_k))\cdot
      \prod_{2\le k\le p,\pm}
\frac{\Gamma(\frac12(q-p)\pm s_k)}{\Gamma(\pm s_k)}
\times\\ &\times
\prod_{2\le k < l\le p,\pm}\frac{\Gamma(\frac12(1+s_k\pm s_l))
          \Gamma(\frac12(1-s_k\pm s_l))}
          {\Gamma(\frac12(s_k\pm s_l))\Gamma(\frac12(-s_k\pm s_l))}
\Phi_s(t)\biggr]
\Biggr|_{s_1=\alpha-\frac12(p+q)+1+2\kappa}ds_2\dots ds_p
\end{align*}
where $\alpha\in \Pi_{\kappa+1}\cap{\cal O}_\delta$.

                         \mal

{\bff 4.\punkt. Calculations.}
 Lemma 4.3.b is an obvious corollary of the last formula.
Nevertheless we present some elements of the calculation,
since it is essential for understanding  Subsection 4.7.

\begin{align*}1)\qquad&\biggl
   [\frac{\Gamma(\frac12(\alpha+\frac12(p+q)-1\pm s_k))}
  {\Gamma(\frac12(s_1\pm s_k))\Gamma(\frac12(-s_1\pm s_k))}
 \biggr] \Bigr|_{s_1=\alpha-\frac12(p+q)+2+2\kappa}
=\\ & =
\frac{\frac12(-\alpha+\frac12(p+q)-1\pm s_k)-\kappa}
     {\Gamma(\frac12(-\alpha+\frac12(p+q)+1\pm s_k))}
\end{align*}
We observe that the  factor $Y(\alpha;s)$
(see (2.2)) is canceled. This factor was the origin
 of singularities in our integral (2.1)--(2.5).

$$2) \qquad \Gamma(\half(1+s_1\pm s_k))
      \Bigr|_{s_1=\alpha-\half(p+q)+1+2\kappa}=
\Gamma(\half(\alpha-\half(p+q)+1+2\kappa\pm s_k))$$
We observe appearance of  factor (4.8) which is very similar to
the factor $Y(\alpha;s)$.
Later it will be an origin of new singularities.

\begin{align*}
3) \qquad&\Bigl[
\frac{\Gamma(\frac12(\alpha-\frac12(p+q)+1+s_1))}
 {\Gamma(s_1)\Gamma(-s_1)}\Bigr]\Bigr|_{s_1=\alpha-\frac12(p+q)+1+2\kappa}=
\\ &=
\frac{(-\alpha+\frac12(p+q)-2\kappa-1)}
{\Gamma(-\alpha+\frac12(p+q)-\kappa)}
\end{align*}

This gives formula (4.6)-(4.12), and  completes the proof
of Lemma 4.3.

\mal

{\bff 4.\punkt. Analytic continuation through the line
  $\Re \alpha<\half(p+q)-1-2\kappa$.}
Lemma 4.3 gives the analytic continuation of $\F_\kappa$ to
${\cal O}_\delta\cap\Pi_{\kappa+1}$. Evidentely the expression for
 the analytic continuation is analytic in the strip
\begin{equation}
-2\kappa-1<\alpha-\half(p+q)+1<-2(\kappa-1)
\end{equation}
and hence we obtain the analytic continuation
of $\F_\kappa$ to the whole strip (4.25).

\mal

{\bff 4.\punkt. Proof of Theorem 2.2.}
The Plancherel
formula (2.1)--(2.5) is correct if $\Re\alpha>\frac12(p+q)-1$.
 We want to construct
 the analytic continuation of its right part to the domain
$\Re\alpha<\frac12(p+q)-1$.  Let us move $\alpha$ to the left side.

 Firstly we pass across the line $\alpha=\frac12(p+q)-1$.
Then we obtain the additional summand
$\F^0(\alpha):=\F_1(\alpha)-\F_0(\alpha)$
given by  formula (4.6)-(4.12) for $\kappa=0$.
This is the summand of the Plancherel formula corresponding
$m=1$, $u_1$=0.

 Let us compare the formula (4.6)--(4.12)
 for $\F^0(\alpha)$  and (2.1)--(2.5).
First, we have in (4.6)--(4.12)  additional factor (4.9). This factor
has singularities but all these singularities are lying in the domain
$\alpha>\frac12(p+q)-1$. The factors (2.2) and (4.8)  are very similar
($\alpha$ is changed to $\alpha+1$). The factors (2.14)--(2.15)
and (4.10)--(4.11) also are very similar. In fact  (4.10)--(4.11)
is the Gindikin--Karpelevich density for $\O(p-1,q-1)$.

Hence, we can construct the analytical continuation of $\F^0(\alpha)$
in the same way as above. The first singularity of $\F^0(\alpha)$ on our way
is the line $\Re\alpha=\frac12(p+q)-2$.
 After passing across the line we obtain one more
summand $\F^{00}(\alpha)$ corresponding $m=2$, $u_1=u_2=0$.

 The line $\Re\alpha=\frac12(p+q)-3$ contains singularities
of the integral $\F(\alpha)$ and also
singularities of $\F^{00}(\alpha)$. Hence, we obtain two additional summands
$\F^{1}(\alpha)$ and $\F^{000}(\alpha)$
 corresponding
$m=1$, $u_1=1$ and $m=3$, $u_1=u_2=u_3=0$
etc. etc. etc.

Formally, we must give complete description of the inductive step but
it literally repeats the arguments of Subsections 4.1--4.6.

\bol

{\bfff 5. Positive definiteness
 of spherical functions.}

\nopagebreak

\bol

\addtocounter{sec}{1}
\setcounter{equation}{0}
\setcounter{fact}{0}
\setcounter{punkt}{0}

We obtained
 the expansion of $\b_\alpha(s)$ in spherical functions
having the form
\begin{equation}
\b_\alpha(z)=\int_{\C^p} \Phi_s(z)\, d\mu_{our}(s)
\end{equation}
where the positive $D_p$-invariant
measure $\mu_{our}(z)$ is described in Theorem 2.2.
Our purpose is to prove positive definiteness of spherical functions
$\Phi_s(z)$ which are contained in the support
of the measure $\mu_{our}$.

By the abstract Plancherel theorem, there exists the unique
 expansion
\begin{equation}
\b_\alpha(z)=\int_{\C^p} \Phi_s(z)\, d\mu_{truth}(s)
\end{equation}
where $\mu_{truth}$ is a positive $D_p$-invariant
 measure on $\C^n$ supported by the space   $\widehat G_{sph}$
of positive definite spherical functions.

Substitute $z=0$ to (5.2). Then $\b_\alpha(0)=1$,
$\Phi_s(0)=1$ and hence
\begin{equation}
\int_{\C^p} d\mu_{truth}=1
\end{equation}

We denote by $\supp \mu_{truth}$ and $\supp \mu_{our}$
the supports of the measures
$ \mu_{truth}$ and $\mu_{our}$.

\mal

{\bff 5.\punkt. Preliminary remarks on the  supports of the measures.}
Consider the bounded polyhedron $Q\subset \R^p$
 described in Theorem 1.14.
Consider the tube $\widetilde Q\subset\C^p$ defined by the condition
$$s\in \widetilde Q \qquad \mbox{iff} \qquad \Re s\in Q$$. Then
\begin{equation}
\supp\mu_{truth}\subset
\widetilde Q ;\qquad \supp\mu_{our}\subset \widetilde Q
\end{equation}
(the first is corollary of Theorem 1.14, the second is corollary of
Theorem 2.2).

Denote by $\R\cup i\R$ the union of the real and imaginary axises in
$\C$.
Then
\begin{equation}
\supp\mu_{our}\subset (\R\cup i\R) \times\dots\times (\R\cup i\R)
 ;\qquad\supp \mu_{truth}\subset (\R\cup i\R) \times\dots\times (\R\cup i\R)
\end{equation}
(the first is the corollary of Theorem 2.2 and the second is
the Corollary of Lemma 1.13).

\mal

{\bff 5.\punkt. Heat kernel.}
Let $\Delta_1$, \dots, $\Delta_p$ be Laplace operators
 (see \cite{Hel}, Section 2.5)
 on the symmetric space
$G/K$. The operator $\Delta_j$ is  some $G$-invariant partial differential
 operator of order $2j$ on $G/K=\B_{p,q}$ with rational
coefficients. The operator $\Delta_1$ is the usual
Laplace-Beltrami operator on $G/K$ (see \cite{Hel}, Section 2.2.4).

 Spherical functions
are joint eigenfunctions of the operators $\Delta_j$
(see \cite{Hel}, Section 4.2). We have  equalities
$$\Delta_j \Phi_s(z)=a_j(s)\Phi_s(z)$$
where $a_j$ are some  polynomials invariant with respect to
the group $D_p$ (see Subsection 1.15) consisting
of permutations and
changing of signs. If $\Phi_s\ne\Phi_{s'}$, then $a_j(s)\ne a_j(s')$
for some $j$.

In particular
$$\Delta_1 \Phi_s(z)=(\lambda+s^2)\Phi_s(z)$$
where $\lambda$ is a constant
and
$$s^2:=s_1^2+\dots+s_p^2$$
By conditions (5.4)--(5.5) {\it the eigenvalues $(\lambda+s^2)$
are  real and they
are uniformly bounded above on the supports of the measures
$\mu_{our}$, $\mu_{truth}$}.

 Consider the Cauchy problem  for the {\it heat equation}
$$(\frac{\partial}{\partial \tau}-\Delta_1)F(z,\tau)=0; \qquad F(z,0)=f(z)$$
on $G/K$.
Let
$R_\tau(z,u)$ be the heat kernel. This means that
the solution of the Cauchy
problem for  the  heat equation  is given by the formula
$$F(z,\tau)=A_\tau f(z):=\int_{G/K} R_\tau(z,u) f(u)\, d\lambda(u)
$$
where $\lambda$ is the  $G$-invariant measure on $G/K$.

\mal

{\scc Lemma 5.\fact.} {\it For each $\tau>0$ and $N$ there exists
a constant $C(\tau,N)$  independent on
$z$, $u$ such that
$$R_\tau(z,u)\le C(\tau,N)
    (1+{\rm dist}(z,u))^{-n}$$
where  $\rm dist(\cdot,\cdot)$ is the distance in $G/K$ associated with
Riemannian metric.}

\mal

{\scc Proof.} Since the heat kernel is $G$-invariant, we can assume
$u=0$. Then
$$R_\tau(z,0)=\int_{s\in i\R}\exp\{\tau(\lambda+s^2)\}
   \Phi_s(z) ds$$
By integral formula for spherical functions (1.38)
$$R_\tau(z,0)=\int_{s\in i\R}\int_{{\bf k}\in K}
\exp\{\tau(\lambda+s^2)\}
                \Psi_s(z^{[{\rm k}]})d{\rm k}\, ds$$
Rapid decreasing of the last expression is more or less obvious.
 \hfill $\boxtimes$

\mal

 Similar estimates are valid for partial derivatives
of  $R_\tau(z,u)$ of any order.

\mal

For spherical functions we have the equality
$$A_\tau\Phi_s(z)=\exp\{ \tau (\lambda+s^2)\}\Phi_s(z)$$

{\scc Lemma 5.\fact.} {\it Let $\mu_{our},\mu_{truth}$
 be the same as above. Then}
\begin{align} a)&
A_\tau \b_\alpha(z)=
          \int_{G/K} \exp\{ \tau (\lambda+s^2)\}\Phi_s(z)\,d\mu_{our}(s)\\
b)&
A_\tau \b_\alpha(z)=
           \int_{G/K} \exp\{ \tau (\lambda+s^2)\}\Phi_s(z)\,d\mu_{truth}(s)
\end{align}

{\scc Proof.} We must prove a possibility to change the order of the
integration. It is sufficient to show absolute convergence of
the integrals
\begin{align}&
\int_{G/K}\int_{\C^p} R_\tau(z,u) \Phi_s(u)\, d\mu_{our}(s)d\lambda(u);
\qquad
&\int_{G/K}\int_{\C^p} R_\tau(z,u) \Phi_s(u) \, d\mu_{truth}(s)d\lambda(u)
\end{align}

 a) For the first integral (5.8)
 we use estimate $|\Phi_s(u)|\le \Phi_{\Re s}(u)$
(see (1.39)).
The measure $\mu_{our}$ is supported
by a finite family of planes $P^{(j)}$ having the form
$$P_j:\qquad \Re s_1=\theta_1^{(j)}, \dots, \Re s_m=\theta_m^{(j)}$$
On each plane $P_j$ the  integrand  is dominated by some expression
$$\sum  R_\tau(z,u) \Phi_{\sigma_j}(u)$$
where $\sigma_j$ are real vectors. The last expression doesn't depend on $s$.
 It remains to notice that
the heat kernel  rapidly decrease
in $u$ for fixed $z$, spherical functions
 $\Phi_{\sigma_j}(u)$ are bounded (see Theorem 1.14),
 and the density of $\mu_{our}(s)$ exponentially
decreases if $|s|\to\infty$.

  b) For the second integral we use the inequality $|\Phi_s(z)|\le 1$
and this completes the proof (see (5.3)). \hfill $\boxtimes$

\mal

{\scc Lemma 5.\fact.}                           {\it
For each polynomial
$r(x_1,\dots,x_p)$
we have                       }
\begin{align*}&
r(\Delta_1,\dots,\Delta_p)A_\tau \b_\alpha(z)=
\int_{G/K}r(a_{1}(s),\dots,a_{p}(s))
 \exp\{ \tau (\lambda+s^2)\}\Phi_s(z)\,d\mu_{our}(s)\\
&r(\Delta_1,\dots,\Delta_p)A_\tau\b_\alpha(z)=
\int_{G/K}r(a_{1}(s),\dots,a_{p}(s))
 \exp\{ \tau (\lambda+s^2)\}\Phi_s(z)\,d\mu_{truth}(s)
\end{align*}

{\scc Proof.}
It is sufficient to prove that all partial derivatives by $z$
of integrals (5.6), (5.7) absolutely converges.
 It is obvious by the following reasons.

1. The integrand rapidly decreases in the variable $s$.

2. For a given $z$ partial derivatives of the heat kernel
by $u$  rapidly decrease.

3. Spherical functions are bounded.\hfill $\boxtimes$

\mal

{\bff 5.\punkt. Proof of positive definiteness.}
Consider  $\sigma\in\supp\mu_{our}$.
Consider the function
$$\eta(s)=\exp\{\tau(\lambda+s^2)\}\sum_{j=1}^p(a_j(s)-a_j(\sigma))^2
    $$
Let $M$ be the maximum of $\eta$
on  $(\R\cup i\R) \times\dots\times (\R\cup i\R)$.
 Then the function
$$\zeta(s)=M-\eta(s)$$
satisfies conditions
$$\zeta(\sigma)=M;\qquad \zeta(s)<M \quad\mbox{if}\quad s\ne w\sigma
\,\,\,\mbox{for all}\,\,\, w\in D_p$$
Consider the sequence  of functions
$$\xi_k(s)=C_k(M-\eta(s))^k\exp\{\tau(\lambda+s^2)\}$$
where $C_k$ is determined by the condition
$\int \xi_k(s)\,d\nu_{our}=1$.
Obviously the sequence $\xi_k(s)$
 converges to distribution $\sum_{w\in D_p}\delta(s-w\sigma)$.

The function $\xi_k(s)$ is a polynomial  expression
$$\xi_k(s)=P_k(a_1(s),\dots,a_p(s),\exp\{\tau(\lambda+s^2)\})$$
Consider the operator
$$\Xi_k:=P_k(\Delta_1,\dots,\Delta_p, A_\tau)$$
By Lemma 5.3 we have
\begin{align}&
\Xi_k\b_\alpha(z)=\int_{G/K}\xi_k(s)\Phi_s(z)d\nu_{our}(s)\\
&\Xi_k\b_\alpha(z)=\int_{G/K}\xi_k(s)\Phi_s(z)d\nu_{truth}(s)
\end{align}
We have $\xi_k(s)\ge 0$. Hence, by Lemma 1.5d and (5.10) the function
 $\Xi_k\b_\alpha(s)$
is positive definite. By (5.9) the sequence
$\Xi_k\b_\alpha(s)$
converges to $\Phi_\sigma(z)$.
Thus, $\Phi_\sigma(z)$ is a point-wise limit of positive definite functions
and hence it is positive definite.

\bol

{\bfff 6. Other series}

\nopagebreak

\bol

\addtocounter{sec}{1}
\setcounter{equation}{0}
\setcounter{fact}{0}
\setcounter{punkt}{0}

{\bff 6.\punkt. Hermitizations.} Below we
 present the  list of hermitizations
(see Subsection 0.1)\footnote%
{Classical part of the list is contained
in H.Jaffee paper \cite{Jaf}. Olshanskii in
\cite{Olshigh}, \cite{Olsdokl} observed that all cases
(including exceptional cases)  can be easily reduced to Nagano
work \cite{Nag}. The list 1-18
is in one-to-one correspondence with the list
of compressive semigroups of symmetric spaces
and with the list of causal symmetric spaces, see
\cite{Olshigh}, \cite{Olsdokl}. The list 1-10 is
in one-to-one correspondence with the list of real
classical categories, see \cite{Ner2}, Addendum A.}.

Each Riemannian noncompact classical
 symmetric space $G/K$ can be realized as a matrix ball%
\footnote{This  observation is present in \cite{Ner}}.
A {\it matrix ball} is a space
of all matrices of a given size over $\K=\R,\C,\H$ with norm $<1$
satisfying (or not satisfying) some symmetry condition.
 The list of matrix balls is given
in the Table 1.

\hfil{\bf Table 1}
\begin{align*}
&\s   \s G/K                  &&\K &&\s\,\,\,\mbox{size}&&\s \mbox{condition}     & &\s \widetilde G/\widetilde K\\
&\s 1. \GL(n,\R)/\O(n,\R)              &&\s\R&&\s n\times n&&\s z=z^t &&\s \Sp(2n,\R)/\U(n)\\
&\s 2. \O(p,q)/\O(p)\times\O(q)        &&\s\R&&\s p\times q&&          &&\s  \U(p,q)/\U(p)\times\U(q)\\
&\s 3. \Sp(2n,\R)/\U(n)                &&\s\C&&\s n\times n&&\s z=z^t  &&\s \Sp(2n,\R)/\U(n)\times\Sp(2n,\R)/\U(n)\\
&\s 4. \GL(n,\C)/\U(n)                 &&\s\C&&\s n\times n&&\s z=z^*  &&\s \U(2n)/\U(n)\times\U(n)\\
&\s 5. \O(n,\C)/\O(n)                  &&\s\R&&\s n\times n&&\s z=-z^t &&\s  \SO^*(2n)/\U(n)\\
&\s 6. \Sp(2n,\C)/\Sp(n)               &&\s\H&&\s n\times n&&\s z=-z^* &&\s  \Sp(4n,\R)/\U(2n)\\
&\s 7. \U(p,q)/\U(p)\times\U(q)        &&\s\H&&\s p\times q&&\s        &&\s [\U(p,q)/\U(p)\times\U(q)]\times [\U(p,q)/\U(p)\times\U(q)]\\
&\s 8. \GL(n,\H)/\Sp(n)                &&\s\H&&\s n\times n&&\s z=z^*  &&\s \SO^*(2n)/\U(n)\\
&\s 9. \Sp(p,q)/\Sp(p)\times\Sp(q)     &&\s\C&&\s p\times q&&\s        &&\s\U(2p,2q)/\U(2p)\times\U(2q)\\
&\s10. \SO^*(2n)/\U(n)                 &&\s\C&&\s n\times n&&\s z=z^t  &&\s \SO^*(2n)/\U(n)\times\SO^*(2n)/\U(n)\\
\end{align*}

The last column contains the hermitization $\widetilde G/\widetilde K$
of $G/K$. The embedding
$G/K\to\widetilde G/\widetilde K$ in all cases is obvious
(we must omit nonholomorphic condition to a matrix $z$, after this we obtain
an hermitian matrix ball).

                   \mal

{\scc Remark.} Some spaces of small dimension are present in the left
column two times (for instance Lobachevskii plane
$\O(2,1)/\O(2)\times\O(1))$. Two associated
hermitizations are different.

                       \mal

Table 2 contains hermitizations related to future tubes

\hfil{\bf Table 2}

\begin{align*}
&\s G/K &&\s  \widetilde G/\widetilde K\\
&11.\s  \SO(2,n)/\O(n)\times\O(2) &&\s  [\SO(2,n)/\O(n)\times\O(2)] \times [\SO(2,n)/\O(n)\times\O(2)]\\
&12.\s  \SO(1,p)\times\O(1,q)  &&\s \SO(2,p+q)
\end{align*}

Exceptional hermitizations%
\footnote{$E{\rm III}$ and $E\rm IV$ are real form of $E_6$,
$E{\rm VII}$
is a real forms of $E_7$, and $F{\rm II}$ is a real form of
$F_4$, see \cite{VO}}
\footnote{Here we observe one of the cases when a phenomenon
existing for classical groups doesn't exist for all exceptional
groups} are given in the Table 3

\hfil{\bf Table 3}
\begin{align*}
&13.\s E{\rm III}/\SO(10)\times \SO(2)        &&\s [E{\rm III}/\SO(10)\times \SO(2)] \times [E{\rm III}/\SO(10)\times \SO(2)]  \\
&14.\s E{\rm VII}/E{\rm III}\times \SO(2)     &&\s  [E{\rm VII}/E{\rm III}\times \SO(2)]  \times [E{\rm VII}/E{\rm III}\times \SO(2)]   \\
&15.\s F{\rm II}/{\rm Spin(9)}                &&\s     E{\rm III}/\SO(10)\times \SO(2)                 \\
&16.\s \Sp(2,2)/\Sp(2)\times\Sp(2)            &&\s   E{\rm III}/\SO(10)\times \SO(2)                 \\
&17.\s \GL(4,\H)/\Sp(2)                       &&\s E{\rm VII}/E{\rm III}\times \SO(2)    \\
&18.\s E{\rm IV}\times\R/F_4                  &&\s E{\rm VII}/E{\rm III}\times \SO(2)
\end{align*}

{\bff 6.\punkt. Kernel representations.}
A {\it kernel representation}%
\footnote{This definition was proposed in \cite{Ner3}} $\rho$ of $G$ is a restriction
of a highest weight representation $\widetilde\rho$
of $\widetilde G$ to $G$. The constructive description
of the scalar-valued kernel
representations of $\O(p,q)$ given in Subsections 1.11-1.12  is valid
for all  series 1-10.

Below we discuss only scalar-valued kernel-representations.

{\bff 6.\punkt. Plancherel formula for large $\alpha$} for the series 1-10
was obtained in \cite{Ner4}%
\footnote{The case of hermitian spaces 3,7,10,11 was considered
by Berezin \cite{Ber2} for $\alpha=\beta$, see notations in Subsection 0.2
(the proof is published in \cite{UU}),
For hermitian case $\alpha\ne\beta$ the Plancherel formula
independently on \cite{Ner4} was obtained by Zhang \cite{Zha}.
Future tube case is simple exercise. The hermitian cases
13, 14 are covered by \cite{UU}; the case 15 is reduced
to one of Gindikin (1964) integrals \cite{Gin1}.
Probably only for the exceptional cases 16, 17, 18 the formula
is not known.
Hence, the possibility to obtain the solution in "general case"
(i.e. 16--18)
is yet preserved.}

  Consider for simplicity nonhermitian case or hermitian
case%
\footnote{For hermitian case $\alpha\ne\beta$ we have the formula
$$C(\alpha)C(\beta)
\prod_{k=1}^p
\Gamma(\frac12(\alpha-h_1+s_k))\Gamma(\frac12(\beta-h_2+s_k))
d{\frak Q}(s)
$$
where $d\frak Q $ is the Shimeno measure \cite{Shi}.  }
$\alpha=\beta$. In all these cases
 the Plancherel formula has the form
\begin{equation}
\prod_{k=1}^p\ch^{-\alpha}t_k=
E(\alpha)\int_{\i\R^p}
\prod_{k\le p,\pm}\Gamma(\frac12(\alpha-h+s_k))\r(s)\Phi_s(t)ds
\end{equation}
where  $m$ is the rank of $G$,
$\r(s)$ is the Gindikin-Karpelevich density,
 $E(\alpha)$ is a meromorphic
factor and  $h$ is a constant. In fact, $h$
 is the last point of square integrability.
This means that
$\int_{G/K}|\b_\alpha(z)|^2dz$ is finite for $\alpha>h$
and infinite for $\alpha=h$.

{\bff 6.\punkt. The analytic continuation of the Plancherel formula.}
 Our arguments from Sections 4--5
don't depend on series and they are valid for all
series 1--10.

In fact considerations of
Section 4 prove
that the following formal procedure gives correct result.
We fix $m=0,1,\dots,p:={\rm rank}\,\,G$ and the collection of
numbers $u_1\le u_2\le\dots u_m$ such that
$\alpha+2u_m+(m-1)\dim \K < h$

Let
\begin{align}
&Q^{m;u}_0(s)=
  \prod_{k\le p,\pm}\Gamma(\frac12(\alpha-h+s_k))\cdot\r(s)\Phi_s(t)\nonumber\\
&Q^{m;u}_{k}(s)
   =\frac{Q_{k-1}^{m;u}(s)}{s_k-\alpha+h-2u_k-k\dim\K}\Biggr|_%
{s_k=\alpha-h+2u_k+k\dim\K}
\end{align}
Then the analytic continuation of (6.1)
has the form
\begin{equation}
\prod_{k=1}^p\ch^{-\alpha}t_k=
\sum\limits_{m;u_1,\dots,u_p} E(\alpha)\frac{(2\pi)^mp!}{(p-m)!}\int_{i\R^{p-m}}
  Q^{m;u}_{m}(s)ds_{m+1}\dots ds_m
\end{equation}
In the cases $\K=\C,\H$ this formula can be considered as a final
result in a closed form.

In the case $\K=\R$ ($G=\O(p,q), \Sp(2n,\R), \GL(n,\R)$)
substitutions (6.3) are impossible without cancellations
and in this case formula (6.3) gives  algorithmic procedure
of calculation of Plancherel measure.

\mal

{\scc Remark. } For the groups $\U(p,q)$ it is easy to obtain
the Plancherel formula (see the explicit
final expression in \cite{Ner3})
 using Berezin--Karpelevich formula for
spherical functions (see \cite{BK},\cite{Hoo}) and Molev
unitarizability results \cite{Molev}; partially this idea was also
realized in \cite{Hil}.

\mal

{\bff 6.\punkt. Kernel representations of compact groups.}
Hermitization procedure also is valid for
 compact Riemannian symmetric spaces.  To obtain the list of hermitizations
we must replace the groups $G$, $\widetilde G$ in Tables 1-3
to their compact forms. For instance the table 1 transforms to the following
table 1'

\pagebreak

\begin{center}
{\bf Table 1'}
\end{center}

         \nopagebreak

\begin{align*}
&\s   \s G/K                   &&\K &&\s\,\,\,\mbox{size}&&\s \mbox{condition}      &&\s \widetilde G/\widetilde K\\
&\s 1'. \U(n)/\O(n,\R)              &&\s\R&&\s n\times n&&\s z=z^t  &&\s \Sp(n)/\U(n)\\
&\s 2'. \O(p+q)/\O(p)\times\O(q)    &&\s\R&&\s p\times q&&          &&\s  \U(p+q)/\U(p)\times\U(q)\\
&\s 3'. \Sp(n)/\U(n)                &&\s\C&&\s n\times n&&\s z=z^t  &&\s \Sp(n)/\U(n)\times\Sp(n)/\U(n)\\
&\s 4'. \U(n)\times\U(n)/\U(n)      &&\s\C&&\s n\times n&&\s z=z^*  &&\s \U(2n)/\U(n)\times\U(n)\\
&\s 5'. \O(n)\times\O(n)/\O(n)      &&\s\R&&\s n\times n&&\s z=-z^t &&\s  \O(2n)/\U(n)\\
&\s 6'. \Sp(n)\times\Sp(n)/\Sp(n)   &&\s\H&&\s n\times n&&\s z=-z^* &&\s  \Sp(2n)/\U(2n)\\
&\s 7'. \U(p+q)/\U(p)\times\U(q)    &&\s\H&&\s p\times q&&\s        &&\s [\U(p+q)/\U(p)\times\U(q)]\times [\U(p+q)/\U(p)\times\U(q)]\\
&\s 8'. \U(2n)/\Sp(n)               &&\s\H&&\s n\times n&&\s z=z^*  &&\s \O(2n)/\U(n)\\
&\s 9'. \Sp(p+q)/\Sp(p)\times\Sp(q) &&\s\C&&\s p\times q&&\s        &&\s\U(2p+2q)/\U(2p)\times\U(2q)\\
&\s10'. \O(2n)/\U(n)                &&\s\C&&\s n\times n&&\s z=z^t  &&\s \O(2n)/\U(n)\times\O(2n)/\U(n)\\
\end{align*}

The construction of the kernel representation of $\O(p+q)$
 given in Subsection 1.28 can be literally translated  to all
 series 1--10. For this purpose we must
replace the space $\M_{p,q}(\R)$ by the space $\M$ of
all matrices over $\K$ (see the second column) having the
size given in the third column and satisfying the condition
in the forth column. The group $G$ acts on  $\M$
 by fractional linear transformations%
\footnote{More details on these models of Riemannian compact
symmetric spaces are contained in \cite{Ner6}}.

Integrals evaluated in \cite{Ner4} easily give
Plancherel formulas for all  kernel-representations
in the cases 1'-10'.

 The case of hermitian symmetric spaces was earlier
considered by Zhang \cite{Zha}.

\mal

{\bff 6.\punkt. Our terminology.} a) {\it The term "hermitization".}
The variant {\it $G/K$ is a real form of} $\widetilde G/\widetilde K$
(see \cite{Jaf}) contradicts to the generally accepted
 usage of the term
"real form", since $\widetilde G/\widetilde K$  is not the
{\it complexification}
$G_\C/K_\C$ of $G/K$.

\mal

b) {\it The term "kernel representation" or "Berezin kernel representation".}
  There is no
common term for this object. One possible
term is {\it "Berezin transform.}
It is not suitable in our situation, since there is no "transform"
in this paper.

Another term  is {\it ``canonical representation''}.
The term ``canonical representation'' can be used in many other senses.
For instance in \cite{VGG} and several successive papers of
Vershik, Gelfand and Graev this term was used for
multiplicative integral and also for some
infinite divisible representations. The complex of phenomena
related to
 infinite divisibility has his own interest
(see, for instance \cite{Ner2}, Chapter 10) but it has
small relation with kernel representations.

There is also the term {\it "Berezin quantization"}. The "quantization"
in Berezin sence is an operation on the space of functions.
In our picture this operation can be defined only for
hermitian case.

\def\itt{\it}

{\sc Moscow State Institute of Electronics and Mathematics
(MIEM)
Bolshoi Triohsviatitelskii per., 3/12,
Moscow -- 109 028,
Russia}

\tt neretin@main.mccme.rssi.ru

\end{document}